\documentclass[oneside,leqno,11pt]{article}
\usepackage{amssymb,amsmath,latexsym}

\def\ga{\mathfrak{a}}

\def\gh{\mathfrak{h}}

\def\gk{\mathfrak{k}}
\def\gl{\mathfrak{l}}

\def\gn{\mathfrak{n}}
\def\go{\mathfrak{o}}
\def\gp{\mathfrak{p}}

\def\gs{\mathfrak{s}}
\def\gt{\mathfrak{t}}
\def\gu{\mathfrak{u}}
\def\gv{\mathfrak{v}}
\def\gw{\mathfrak{w}}

\def\gz{\mathfrak{z}}

\def\ul{\mathbf{l}}
\def\um{\mathbf{m}}

\def\uw{\mathbf{w}}

\def\ggg{> \hskip -5 pt >}

\def\Aut{{\rm Aut}}

\def\Ad{{\rm Ad}}

\def\Ind{{\rm Ind\,}}

\def\Im{{\rm Im\,}}
\def\Re{{\rm Re\,}}
\def\Pf{{\rm Pf}}

\def\Skew{{\rm Skew\,}}  % antisymmetric matrices
\def\C{\mathbb{C}}

\def\F{\mathbb{F}}

\def\H{\mathbb{H}}

\def\O{\mathbb{O}}
\def\P{\mathbb{P}}

\def\R{\mathbb{R}}

\def\cA{\mathcal{A}}

\def\cE{\mathcal{E}}
\def\cF{\mathcal{F}}

\def\cH{\mathcal{H}}

\def\cL{\mathcal{L}}
\def\cM{\mathcal{M}}

\def\cO{\mathcal{O}}

\def\cU{\mathcal{U}}

\newtheorem{theorem}[equation]{Theorem}

\newtheorem{lemma}[equation]{Lemma}

\newtheorem{corollary}[equation]{Corollary}

\newtheorem{proposition}[equation]{Proposition}

\newtheorem{definition}[equation]{Definition}

\newtheorem{remark}[equation]{Remark}

\title{Infinite Dimensional Multiplicity Free Spaces II: \\
Limits of Commutative Nilmanifolds}

\begin{document}
\author{Joseph A. Wolf} 
\date{24 January, 2008}

\maketitle

\begin{abstract}
We study direct limits $(G,K) = \varinjlim\, (G_n,K_n)$ of Gelfand 
pairs of the form $G_n = N_n\rtimes K_n$ with $N_n$ nilpotent,
in other words pairs $(G_n,K_n)$ for which $G_n/K_n$ is a commutative
nilmanifold.  First, we extend the criterion of \cite{W3} for a direct 
limit representation to be multiplicity free.  Then we study direct 
limits $G/K = \varinjlim \, G_n/K_n$ of commutative nilmanifolds and
look to see when the regular representation of $G = \varinjlim G_n$ on 
an appropriate Hilbert space $\varinjlim L^2(G_n/K_n)$ is multiplicity free.  
One knows that the $N_n$ are commutative or $2$--step nilpotent.  In many
cases where the derived algebras $[\gn_n,\gn_n]$ are of bounded dimension 
we construct $G_n$--equivariant isometric maps 
$\zeta_n : L^2(G_n/K_n) \to L^2(G_{n+1}/K_{n+1})$
and prove that the left regular representation of $G$ on the Hilbert space
$L^2(G/K) := \varinjlim \{L^2(G_n/K_n),\zeta_n\}$ is a multiplicity free
direct integral of irreducible unitary representations.  The direct
integral and its irreducible constituents are described explicitly.
One constituent of our argument is an extension of the classical
Peter--Weyl Theorem to parabolic direct limits of compact groups.
\end{abstract}

\section{Introduction} \label{sec1}
\setcounter{equation}{0}

Gelfand pairs $(G,K)$, and the corresponding ``commutative'' homogeneous
spaces $G/K$, form a natural extension of the class of riemannian symmetric
spaces.  
Let $G$ be a locally compact topological group, $K$ a compact subgroup,
and $M = G/K$.
Then the following conditions are equivalent; see \cite[Theorem 9.8.1]{W2}.
\medskip

\noindent\phantom{XX}
1. $(G,K)$ is a Gelfand pair, i.e. $L^1(K\backslash G/K)$ is
commutative under convolution.
\hfill\newline\phantom{XX}
2. If $g,g' \in G$ then $\mu_{_{KgK}} * \mu_{_{Kg'K}} =
 \mu_{_{Kg'K}} * \mu_{_{KgK}}$
{\rm (convolution of measures on $K\backslash G/K$).}
\hfill\newline\phantom{XX}
3. $C_c(K\backslash G/K)$ is commutative under convolution.
\hfill\newline\phantom{XX}
4. The measure algebra $\cM(K\backslash G/K)$ is commutative.
\hfill\newline\phantom{XX}
5. The representation of $G$ on $L^2(M)$ is multiplicity free.
\medskip

\noindent If $G$ is a connected Lie group one can also add
\medskip

\noindent\phantom{XX}
6. The algebra of $G$--invariant differential operators on $M$ is commutative.
\medskip

\noindent
Conditions 1, 2, 3 and 4 depend on compactness of $K$ so that integration 
on $M$ and $K\backslash G/K$ corresponds to integration on $G$.  Condition 
5 makes sense as long as $K$ is unimodular in $G$, and condition 6 remains 
meaningful (and useful) whenever $G$ is a connected Lie group.
\medskip

In this note we look at some cases where $G$ and $K$ are not locally
compact, in fact are infinite dimensional, and show in those cases that
the multiplicity--free condition 5 is satisfied.  The main results
are Theorems \ref{heis-case} and \ref{big-hnt}.  We first discuss a
multiplicity free criterion which extends that of \cite{W3}.  We used 
that criterion in \cite{W3} for a number of direct systems $\{(G_n,K_n)\}$ 
where the $G_n$ are compact Lie groups and the $K_n$ are closed subgroups.
Here we use its extension for direct systems $\{(G_n,K_n)\}$ in which a
(closed connected) nilpotent subgroup $N_n$ of $G_n$ acts transitively
on $G_n/K_n$.  Then (see \cite[Chapter 13]{W2}) $N_n$ is the nilradical
of $G_n$ and $G_n$ is the semidirect product group $N_n\rtimes K_n$.  Further,
the group $N_n$ is commutative or $2$--step nilpotent so its structure and
representation theory can be clearly understood in terms of the familiar
Heisenberg groups. 
\medskip

Section \ref{sec3} pins down the structure of direct limit of Heisenberg
group, and Section \ref{sec4} provides an infinite dimensional analog of
the Peter--Weyl Theorem that turns out to be a necessary tool for our 
multiplicity free arguments.
\medskip

In Section \ref{sec5} we look closely at the cases where $N_n$ is the
$(2n+1)$--dimensional Heisenberg group $H_n = \Im\C + \C^n$ and $K_n$
(viewed as a subgroup of $U(n)$) acts irreducibly on $\C^n$.  The
classification of these cases is well known; see Table \ref{kac-table}.
It leads to $16$ direct systems, listed in Table \ref{jaw-table}.  In
Theorem \ref{heis-case} we see that, for 
each of these systems, the limit space $L^2(G/K)$ is well 
defined and the action of $G$ on $L^2(G/K)$ is multiplicity free.
The methods here are explicit and they guide the arguments in the
more general cases.
\medskip

Our basic arguments go in two steps.  First we examine the limit
$L^2(N) := \varinjlim L^2(N_n)$ as an $(N\times N)$--module, and then
we look at the right action of $K$ on $L^2(N)$ in order to analyze
$L^2(G/K) := \varinjlim L^2(G_n/K_n)$ as a left $G$--module.  The first 
step is carried out in Section \ref{sec6}.  Section \ref{sec7} works out
some preliminaries for understanding the right action of $K$.  Then in
Sections \ref{sec8} and \ref{sec9} specify the precise requirements
for the multiplicity free condition.  See Theorems \ref{hnt3} and 
\ref{iso-gelfand}.  The latter reduces some considerations to the case
where $N_n$ is a Heisenberg group.
\medskip

Tables \ref{ind-vin-table} and \ref{ind-vin-table-ipms} list the key cases 
of direct systems of nilpotent commutative pairs not based directly on
Heisenberg groups.  The corresponding multiplicity free results are
Theorems \ref{big-hnt} and \ref{big-hhhhnt}.
\medskip

Finally, Appendix A is a small discussion of formal degree for induced
representations, and Appendix B indicates how some of our results can be
made more explicit using branching rules for representations of the $K_n$.
\medskip

Our arguments require the direct system $\{(G_n,K_n)\}$ to have the
property that the $\dim [\gn_n,\gn_n]$ have an upper bound.  It seems
probable that this condition is not necessary, but that will require a
new idea.

\section{Direct Limit Groups and Representations}\label{sec2}
\setcounter{equation}{0}

We consider direct limit groups $G = \varinjlim G_n$ and direct limit
representations $\pi = \varinjlim \pi_n$ of them.  This means that
$\pi_n$ is a representation of $G_n$ on a vector space $V_n$, that
the $V_n$ form a direct system, and that $\pi$ is the representation 
of $G$ on $V =\varinjlim V_n$ given by $\pi(g)v = \pi_n(g_n)v_n$ whenever
$n$ is sufficiently large that $V_n \hookrightarrow V$ and $G_n \hookrightarrow
G$ send $v_n$ to $v$ and $g_n$ to $g$.  The formal definition amounts to
saying that $\pi$ is well defined.
\medskip

Here we will only consider unitary representations.  Thus the $V_n$ all
will be Hilbert spaces and the inclusions $V_n \hookrightarrow V_{n+1}$
will preserve norms (and inner products).

It is clear that a direct limit of irreducible representations is
irreducible, but there are irreducible representations of direct limit
groups that cannot be formulated as direct limits of irreducible
finite dimensional representations.  This is a combinatoric matter and
is discussed extensively in \cite{DPW}.  Dealing with this matter is
the crux of the problem of proving multiplicity--free properties.  
The following definition is closely related to the relevant combinatorics 
but applies to a somewhat simpler situation.

\begin{definition}\label{limit--aligned}{\rm
We say that a unitary representation $\pi$ of $G = \varinjlim G_n$ is
{\bf limit--aligned} if it is a direct limit $\varinjlim \pi_n$
of unitary representations in such a way that (i) each group $G_n$ is
of type I, (ii) $\pi_n$ is a continuous direct sum 
$\int \zeta_n\ d\nu_n(\zeta_n)$ of mutually disjoint primary representations, 
and (iii) in the corresponding inclusions 
$V_{\pi_n} = \int V_{\zeta_n}\ d\nu_n(\zeta_n) \hookrightarrow 
\int V_{\zeta_{n+1}}\ d\nu_{n+1}(\zeta_{n+1}) = V_{\pi_{n+1}}$ of
representation spaces map $\nu_n$--almost-every primary integrand 
$V_{\zeta_n}$ of $V_{\pi_n}$
into a primary integrand $V_{\zeta_{n+1}}$ of $V_{\pi_{n+1}}$.
}
\end{definition}

\begin{theorem}\label{reduction}
Let $\pi = \varinjlim \pi_n$ be a limit--aligned unitary representation
of $G = \varinjlim G_n$.  Suppose that the $\pi_n$ are multiplicity free.
Then $\pi$ is multiplicity free.  In other words the commuting algebra 
of $\pi$ is commutative.
\end{theorem}

\noindent {\bf Proof.}  
Let $V = \varinjlim V_{\pi_n}$ be the representation spaces.  Consider the
primary decompositions $V_{\pi_n} = \int V_{\zeta_n}\ d\nu_n(\zeta_n)$.
Here $\nu_n$ is a non--negative measure on the analytic Borel space
$\widehat{G_n}$, the $\zeta_n$ are mutually inequivalent irreducible 
unitary representations of $G_n$, and $\pi_n = \int \zeta_n\ d\nu_n(\zeta_n)$.
\medskip

Since $\pi$ is limit--aligned we have maps $b_n : (\widehat{G_n},\nu_n)
\to (\widehat{G_{n+1}},\nu_{n+1})$ of measure spaces (modulo sets of
measure zero) such that, for $\nu_n$--almost all $\zeta_n \in \widehat{G_n}$,
$V_{\pi_n} \hookrightarrow V_{\pi_{n+1}}$ maps $V_{\zeta_n}$ into
$V_{b_n(\zeta_n)}$.  Thus we have direct systems
$$
V_{\zeta_n} \to V_{b_n(\zeta_n)} \to V_{b_{n+1}(b_n(\zeta_n))}
	\to V_{b_{n+2}(b_{n+1}(b_n(\zeta_n)))} 
	\to V_{b_{n+3}(b_{n+2}(b_{n+1}(b_n(\zeta_n))))} \to \dots
$$
of irreducible unitary representation spaces for $G = \varinjlim G_n$.
Let $V_\zeta$ denote the corresponding direct limit Hilbert space and
$\zeta$ the (necessarily irreducible) direct limit unitary representation
of $G$ on $V_\zeta$.
\medskip

Write $\widehat{G_n}'$ for the support of $\nu_n$ in the hull--kernel 
topology.  From the considerations just above we see that 
$b_n(\widehat{G_n}') \subset \widehat{G_{n+1}}'$.  That defines a 
space $\widehat{G}' = \varinjlim \widehat{G_n}'$, a measure class
$\nu = \varinjlim \nu_n$ on $\widehat{G}'$, and a decomposition
$V = \int V_\zeta\ d\nu(\zeta)$.  Since the $\zeta$ are irreducible,
the closed $\pi(G)$--invariant subspaces of $V$ are just the
$V_S = \int_S V_\zeta\ d\nu(\zeta)$ where $S$ is a measurable subset
of $\widehat{G}'$.  Thus any two projections in the commuting algebra
of $\pi$ commute with each other.  Since the commuting algebra is
a $W^*$ algebra, thus generated by projections, it is commutative.
We have proved that $\pi$ is multiplicity free.  \hfill $\square$
\medskip

\section{Direct Limits of Heisenberg Groups}\label{sec3}
\setcounter{equation}{0}

In this section we work out the structure and properties of the 
($2$--sided) regular representation of the infinite Heisenberg groups.
That is the foundation for study of the multiplicity free
property for direct limits of commutative nilmanifolds.
\medskip

Recall that the (ordinary) Heisenberg group is the group $H_n = \Im\C + \C^n$
with composition $(z,w)(z',w') = (z + z' + \Im(w\cdot w'), w+w')$.  Here
$w\cdot w'$ refers to the standard positive definite hermitian inner product
on $\C^n$ and $\Im v = \tfrac{1}{2}(v-\overline{v})$ is the imaginary
{\em component} of a complex number $v$.  Thus $\Im\C$ is both the center
and the derived group and $H_n/\Im\C$ is a vector group.
\medskip

For $t$ nonzero and real, 
$\pi_{n,t}$ denotes the irreducible unitary representation of $H_n$ with 
central character $e^t: (z,0) \mapsto e^{tz}$ (we use the fact that 
$z \in \Im\C$ is pure imaginary).  Then $\pi_{n,t}$ is square integrable 
modulo the center of $H_n$ in the sense that its coefficients 
$f_{u,v}(h) = \langle u , \pi_{n,t}(h)v\rangle$
satisfy $|f| \in L^2(H_n/\Im\C)$.  For the appropriate normalization of 
Haar measure $\pi_{n,t}$ has formal degree $|t|^n$ in the sense of the
orthogonality relation $\langle f_{u,v},f_{u',v'}\rangle_{L^2(H_n/\Im\C)}
= |t|^{-n}\langle u,u'\rangle\, \overline{\langle v,v'\rangle}$.
\medskip

The representation space of $\pi_{n,t}$ is the Fock space
$$
\cH_{n,t} = \left \{f : \C^n \to \C \text{ holomorphic} \ \left | \
\int |f(w)|^2\exp(-|t||w|^2) d\lambda(w) < \infty\right \}\right .
$$
where $\lambda$ is Lebesgue measure.  The representation is
$$
[\pi_t(z,v)f](w) = e^{tz \pm t\Im(w- v/2)\cdot v} f(w-v)
$$ 
where $\pm $ is the sign of $t/|t|$.
\medskip

For each multi--index
$\um = (m_1, \dots , m_n)$, $m_i \geqq 0$,  we have the the monomial
$w^{\um} = w_1^{m_1}\dots w_n^{m_n}$.
One computes $\int w^{\um} \overline{w^{\um'}}\exp(-|t||w|^2) d\lambda(w)$ to
see that it is equal to $0$ for $\um \ne \um'$, and if $\um = \um$ it is 
equal to $c\,\um !$ for a positive constant $c$ independent of $\um$.
Here $\um !$ means $\prod (m_i!)$.  Thus we can (and do) normalize the inner
product on $\cH_t$ so that the $w[\um] := w^{\um}/\sqrt{\um !}$ form a complete
orthonormal set in $\cH_{n,t}$.  The corresponding space of matrix coefficients
is $\cE_{n,t} = \cH_{n,t} \widehat{\otimes} \cH_{n,t}^*$.  It is spanned by
the $f_{\ul,\um;t}: g \mapsto \langle w[\ul],\pi_{n,t}(g)z[\uw]\rangle$.
These coefficients belong to the Hilbert space
$$
L^2(H_n/\Im\C;e^t) = \{f:H_n \to \C \mid |f| \in L^2(H_n/\Im\C)
	\text{ and } f(z,w) = e^{-tz}f(0,w)\}
$$
with inner product $\langle f,f'\rangle = \int_{\C^n} f(z,w)\overline{f'(z,w)}
d\lambda(w)$.
\medskip

Since $|t|^n$ is the formal degree of $\pi_{n,t}$ the orthogonality
relations say that the inner product in $\cE_{n,t}$ is given by
$\langle f_{\ul,\um;t}, f_{\ul',\um';t}\rangle = |t|^{-n}$ if $\ul = \ul'$
and $\um = \um'$, $0$ otherwise.  Now the $|t|^{n/2}f_{\ul,\um;t}$
form a complete orthonormal set in $\cE_{n,t}$, and $\cE_{n,t}$ consists of
the functions $\Phi_{n,t,\varphi}$ given by 
\begin{equation}
\Phi_{n,t,\varphi}(h) =
\sum_{\ul,\um} \varphi_{\ul,\um}(t) |t|^{n/2} f_{\ul,\um;t}(h)
\end{equation} 
where the numbers $\varphi_{\ul,\um}(t)$ satisfy
$\sum_{\ul,\um} |\varphi_{\ul,\um}(t)|^2 < \infty$.
\medskip

The Hilbert space $L^2(H_n)$ is the direct integral 
$\int_{-\infty}^\infty \cE_{n,t} |t|^n\, dt$ based on the complete
orthonormal sets $\{|t|^{n/2}f_{\ul,\um;t}\}$ in the $\cE_{n,t}$.  
Thus it consists of all functions $\Psi_{n,\varphi}$ given by
\begin{equation}\label{a1}
\Psi_{n,\varphi}(h) = \int_{-\infty}^\infty \Phi_{n,t,\varphi}(h)|t|^n\, dt
= \int_{-\infty}^\infty \left ( {\sum}_{\ul,\um}
\varphi_{\ul,\um}(t) |t|^{n/2}f_{\ul,\um;t}(h) \right ) |t|^n\, dt
\end{equation} 
such that the functions $\varphi_{\ul,\um}: \R \to \C$ are measurable,
$\sum_{\ul,\um} |\varphi_{\ul,\um}(t)|^2 < \infty$ for almost all $t$, and
$\sum_{\ul,\um} |\varphi_{\ul,\um}(t)|^2 \in L^1(\R,|t|^ndt)$.
Note that 
\begin{equation}\label{a2}
\begin{aligned}
||\Psi_{n,\varphi}||^2_{L^2(H_n)} &= \int_{-\infty}^\infty
	||\Phi_{n,t,\varphi}||^2_{\cE_{n,t}}\, |t|^n\, dt\\
	&= \int_{-\infty}^\infty 
  	  \left ({\sum}_{\ul,\um} |\varphi_{\ul,\um}(t)|^2\right ) |t|^n\, dt
= {\sum}_{\ul,\um} ||\varphi_{\ul,\um}||^2_{L^2(\R,|t|^{n/2}dt)}
\end{aligned}
\end{equation}

The left/right representation of $H_n\times H_n$ on $\cE_{n,t}$
is the exterior tensor product $\pi_{n,t} \boxtimes \pi_{n,t}^*$.  It is
irreducible.  The corresponding representation of $H_n\times H_n$ on
$L^2(H_n)$ is the left/right regular representation 
$\Pi_n:= \int_{-\infty}^\infty (\pi_{n,t} \boxtimes \pi_{n,t}^*)|t|^ndt$.
(The factor $|t|^n$ is not relevant to the equivalence class of the
representation, but it is crucial to expansion of functions.)  
\begin{lemma}\label{regrephn}
The left/right regular representation $\Pi_n$ of $H_n\times H_n$ on
$L^2(H_n)$ is multiplicity free.
\end{lemma}

\noindent {\bf Proof.} Any invariant subspace must be invariant under the
center $\Im \C \times \Im \C$ of $H_n\times H_n$, so it is of the form
$\int_S \cE_{n,t}\, |t|^n\, dt$ where $S$ is a measurable subset of
$\R$.  Now any two projections in the commuting algebra of $\Pi_n$ must
commute with each other.  The commuting algebra is a von Neumann algebra,
so the projections generate a dense subalgebra.  Thus the
commuting algebra is commutative.  \hfill $\square$
\medskip

Suppose $m \geqq n$.  We view an $n$--tuple $\uw$ as an $m$--tuple by
appending $m-n$ zeroes; then $w^{\um}$ and $w[\um]$ have the same meaning
as functions on $\C^n$ and on $\C^m$.  Thus also the coefficient function
$f_{\ul,\um;t}: H_n \to \C$ is the restriction of $f_{\ul,\um;t}: H_m \to \C$. 
From (\ref{a1}) and (\ref{a2}) we see that the map
\begin{equation}\label{a3}
\zeta'_{m,n}: \cE_{n,t} \to \cE_{m,t} \text{ defined by }
\zeta'_{m,n}(|t|^{n/2}f_{\ul,\um;t}) = |t|^{m/2}f_{\ul,\um;t} \text{ and }
\zeta'_{m,n}(\Psi_{n,\varphi}) = \Psi_{m,|t|^{(n-m)/2}\varphi}
\end{equation}
is an isometric $(H_n\times H_n)$--equivariant injection of 
$\cE_{n,t}$ into $\cE_{m,t}$.  Specifically, it maps a complete
orthonormal set in $\cE_{n,t}$ to an orthonormal set in $\cE_{m,t}$.
 
\begin{remark}\label{res-heis}{\rm
The adjoint of the isometric injection $\zeta'_{m,n}: \Psi_{n,\varphi} \mapsto
\Psi_{m,|t|^{(n-m)/2}\varphi}$ of $L^2(H_n)$ into $L^2(H_m)$ is 
orthogonal projection of $L^2(H_m)$ to the image of
the injection.  It is the scalar
multiple (by $|t|^{(m-n)/2}$) of restriction of functions on each
direct integrand $\cE_{m,t}$ of $L^2(H_m)$, but clearly the scalar varies 
with $m$, $n$ and $t$.}  \hfill $\diamondsuit$
\end{remark}

Now let's consider convergence.  Let $M_n$ denote the set of 
multi--indices $\um$ of length $n$.  In order that $\Phi_{n,t} \in \cE_{n,t}$
we needed that $\varphi(t) \in L^2(M_n)$, using counting measure.
That gives the function $||\varphi||^2(t)$.  Now for 
$\Psi_{n,\varphi}$ to be in $L^2(H_n)$
as $n$ grows, we need $||\varphi||^2 \in L^1(\R,|t|^ndt)$ as $n$ grows.
These conditions are satisfied whenever only finitely many of the
$\varphi_{\ul,\um}$ are not identically zero, and the nonzero ones
are Schwartz class functions of $t$.  Thus we have a dense subset of
$L^2(H_n)$ that maps in a norm--preserving way into $L^2(H_m)$,
whenever $m \geqq n$, and extends by continuity to give a well defined 
unitary injection $L^2(H_n) \to L^2(H_m)$.  
\medskip

All the ingredients in the
construction of this injection are $(H_n\times H_n)$--equivariant, so the
unitary injection $L^2(H_n) \to L^2(H_m)$ is equivariant for the
left/right regular representation $\Pi_n$ of $H_n\times H_n$.  

\begin{theorem}\label{heis-inject}
There is a strict direct system $\{L^2(H_n),\zeta'_{m,n}\}$ of $L^2$ 
spaces of the
Heisenberg groups, whose maps $\zeta'_{m,n}: L^2(H_n) \to L^2(H_m)$ are
$(H_n\times H_n)$--equivariant unitary injections.  Let $\Pi_n$ denote the 
left/right regular representation of $H_n\times H_n$ on $L^2(H_n)$
and let $H$ denote the infinite dimensional Heisenberg group 
$\varinjlim H_n$.  Then we have a well defined Hilbert 
space $L^2(H): = \varinjlim \{L^2(H_n), \zeta'_{m,n}\}$ and a natural unitary 
representation $\Pi = \varinjlim \Pi_n$ of $H \times H$ on
$L^2(H)$.  Further, that representation $\Pi$ is multiplicity--free.
\end{theorem}

\noindent {\bf Proof.} All the assertions except the multiplicity--free 
assertion have just been proved.  Now it remains only to prove that 
$\Pi = \varinjlim \Pi_n$ is limit--aligned, for then Lemma \ref{regrephn} 
and Theorem \ref{reduction} complete the proof.  That alignment is
immediate because the unitary injections $L^2(H_n) \to L^2(H_m)$
are equivariant for the action of the center of $H_n\times H_n$, and that
action specifies the direct integrands $\cE_{n,t}$ and $\cE_{m,t}$ 
within $L^2(H_n)$ and $L^2(H_m)$.
\hfill $\square$

\section{The Peter--Weyl Theorem for Direct Limits of Compact Groups}
\label{sec4}
\setcounter{equation}{0}

We will apply Theorem \ref{heis-inject} to the study of direct limits
$\varinjlim \{(H_n\rtimes K_n,K_n)\}$.  There $K_n$ is a compact connected 
group of automorphisms of the Heisenberg group $H_n$.  In order to do that
we need to extend Theorem \ref{heis-inject} from the $H_n$ to the 
semidirect product groups $H_n\rtimes K_n$, and for that we must first
prove the analog of Theorem \ref{heis-inject} for the direct systems
$\{K_n\}$.  This analog, the Peter--Weyl Theorem for $\{K_n\}$, is the
subject of this section.  As we will see in Section \ref{sec5} we need
only study the restricted class of direct systems $\{K_n\}$ given as follows.
\medskip

\begin{definition}\label{paralim} {\rm
Let $\{K_n\}$ be a strict direct system of compact connected Lie groups,
$\{(K_n)_{_\C}\}$ the direct system of their complexifications.
Suppose that, for each $n$,
\begin{equation}
\begin{aligned}
&\text{\,the semisimple part } \ [(\gk_n)_{_\C}, (\gk_n)_{_\C}] \
\text{ of the reductive algebra } \ (\gk_n)_{_\C} \\
&\text{ is the semisimple component of a parabolic
subalgebra of } (\gk_{n+1})_{_\C}.
\end{aligned}
\end{equation}
Then we say that the direct systems $\{K_n\}$ and $\{(K_n)_{_\C}\}$
are {\bf parabolic} and that $\varinjlim K_n$ and $\varinjlim (K_n)_{_\C}$
are {\bf parabolic direct limits}.  This is a small variation on
the definitions of parabolic and weakly parabolic direct limits in
\cite{W1}. }\hfill $\diamondsuit$
\end{definition}

Now let $\{K_n\}$ be a strict direct system of compact connected Lie groups
that is parabolic.  We recursively construct Cartan subalgebras
$\gt_n \subset \gk_n$ with $\gt_1 \subset \gt_2 \subset \dots \subset 
\gt_n \subset \gt_{n+1} \subset \dots$ and (using the parabolic
property) simple root systems
$\Psi_n = \Psi((\gk_n)_{_\C}, (\gt_n)_{_\C})$ such that 
each simple root for $(\gk_n)_{_\C}$ is the restriction of exactly one
simple root for $(\gk_{n+1})_{_\C}$.  Then we may assume that
$\Psi_n = \{\psi_{n,1}, \dots ,\psi_{n,p(n)}\}$ in such a way that each
$\psi_{n,j}$ is the $(\gt_n)_{_\C}$--restriction of $\psi_{n+1,j}$
and of no other element of $\Psi_{n+1}$.  The corresponding sets 
$\Xi_n = \{\xi_{n,1}, \dots , \xi_{n,p(n)}\}$ of of fundamental highest
weights satisfy: $\xi_{n+1,j}$ is the unique element of $\Xi_{n+1}$ whose
$(\gt_n)_{_\C}$--restriction is $\xi_{n,j}$, for $1 \leqq j \leqq p(n)$.
\medskip

Write $\kappa_{n,\lambda}$ for the irreducible representation of $K_n$
with highest weight $\lambda_n = \sum_1^{p(n)} \ell_{n,j}\xi_{n,j}$
and $\cF_{n,\lambda}$ for its representation space.  Choose highest 
weight unit vectors $v_{n,\lambda} \in \cF_{n,\lambda}$.  Note that
the Lie algebra inclusion $\gk_n \hookrightarrow \gk_{n+1}$ defines 
an inclusion $\cU(\gk_n) \hookrightarrow  \cU(\gk_{n+1})$ of enveloping
algebras.  Now $\cA(v_{n,\lambda}) \mapsto \cA(v_{n+1,\lambda})$,
$\cA \in \cU(\gk_n)$, defines a $K_n$--equivariant injection
$\cF_{n,\lambda} \hookrightarrow \cF_{n+1,\lambda}$, and that injection
is unitary because it sends the unit vector $v_{n,\lambda}$ to the
unit vector $v_{n+1,\lambda}$.  Thus we can simply regard
$\cF_{n,\lambda}$ as the subspace of $\cF_{n+1,\lambda}$ generated by the
highest weight unit vector $v_{n,\lambda} = v_{n+1,\lambda}$.
\medskip

The matrix coefficients of $\kappa_{n,\lambda}$ are the 
$f_{x,y}(k) = \langle x, \kappa_{n,\lambda}(k)y\rangle$, and the
Schur Orthogonality Relations say that their $L^2(K_n)$ inner product
is $\langle f_{x,y},f_{x',y'}\rangle = \deg(\kappa_{n,\lambda})^{-1}
\langle x,x'\rangle \overline{\langle y,y'\rangle}$.  Write
$\cL_{n,\lambda}$ for that space of matrix coefficients.  For $m\geqq n$
let 
$$
\zeta''_{m,n} : f_{x,y} \text{ (as a function on } K_n) \mapsto 
(\deg(\kappa_{m,\lambda})/\deg(\kappa_{n,\lambda}))^{1/2}f_{x,y}
\text{ (as a function on } K_m).
$$  Then 
$$
\zeta''_{m,n}: \cL_{n,\lambda} \to \cL_{m,\lambda} \text{ is a }
(K_n\times K_n)\text{--equivariant isometric map.} 
$$
The Peter--Weyl Theorem $L^2(K_n) = \sum_\lambda \cL_{n,\lambda}$ assembles the
$\cL_{n,\lambda} \hookrightarrow \cL_{n+1,\lambda}$ into a
$(K_n\times K_n)$--equivariant isometric injection 
$\zeta''_{n,m}: L^2(K_n) \to
L^2(K_m)$.  We now come to the following result, corresponding to 
Theorem \ref{heis-inject}.

\begin{theorem}\label{lim-pw}
{\rm (Peter--Weyl Theorem for parabolic direct limit groups.)\ }
Let $\{K_n\}$ be a strict direct system of compact connected Lie
groups.  Suppose that $\{K_n\}$ is parabolic.  Then
there is a strict direct system $\{L^2(K_n), \zeta''_{m,n}\}$ of $L^2$ spaces, 
whose maps $\zeta''_{m,n}: L^2(K_n) \to L^2(K_m)$ are
$(K_n\times K_n)$--equivariant unitary injections.  Let $\Gamma_n$ denote the
left/right regular representation of $K_n\times K_n$ on $L^2(K_n)$
and let $K = \varinjlim K_n$.  Then we have a well defined Hilbert
space $L^2(K): = \varinjlim \{L^2(K_n), \zeta''_{m,n}\}$ and a natural unitary
representation $\Gamma = \varinjlim \Gamma_n$ of $K \times K$ on
$L^2(K)$.  

The space $L^2(K)$ is the Hilbert space orthogonal direct sum of
its subspaces $\cL_\lambda := \varinjlim \cL_{n,\lambda}$ and the action
of $K \times K$ on $\cL_\lambda$ is irreducible with highest weight $\lambda$.  
In particular the left/right regular representation $\Gamma$ is 
multiplicity--free.
\end{theorem}

\noindent {\bf Proof.} As in the Heisenberg group setting, it remains 
only to prove that $\Gamma = \varinjlim \Gamma_n$ is limit--aligned, for then
Lemma \ref{regrephn} and Theorem \ref{reduction} complete the proof.  
Evidently $L^2(K)$ contains the mutually orthogonal spaces
$\cL_\lambda := \varinjlim \cL_{n,\lambda}$.  If $f \in L^2(K)$
is orthogonal to all of them, we interpret $f$ as a function on $K$,
and then $f|_{K_n} = 0$ for all $n$, so $f = 0$.  Thus
$L^2(K)$ is the (Hilbert space closure of the) orthogonal direct 
sum of the $\cL_\lambda := \varinjlim \cL_{n,\lambda}$. As
$\lambda$ determines the direct summand  $\cL_\lambda$, now $\Gamma$
is limit--aligned.
\hfill $\square$

\begin{remark}\label{res-compact}{\rm 
Just as we noted in Remark \ref{res-heis}, the adjoint to the injection
$\zeta''_{m,n}: L^2(K_n) \to L^2(K_m)$ is orthogonal projection
to the image of that injection, and on each $(K_m \times K_m)$--irreducible
summand of $L^2(K_m)$ it is a scalar multiple of restriction of functions.}
\hfill $\diamondsuit$
\end{remark}

\section{Commutative Spaces for Heisenberg Groups}\label{sec5}
\setcounter{equation}{0}

In this section we put together the results of Sections \ref{sec3} and
\ref{sec4} to study direct systems $\{(G_n,K_n)\}$ of Gelfand pairs 
where $G_n$ is the semidirect product $H_n\rtimes K_n$ of a Heisenberg
group with a compact subgroup $K_n \subset \Aut(H_n)$.  Then $K_n$ 
is a closed subgroup of the maximal compact subgroup $U(n) \in \Aut(H_n)$.
The concrete results in this section will require that $K_n$ be
connected and that its action on $\C^n$ be irreducible.
\medskip

The classification goes as follows for the cases where $K_n$ is 
connected and is irreducible 
on $\C^n$.  Carcano's theorem (\cite{C}; or see 
\cite[Theorem 4.6]{BJR} or \cite[Theorem 13.2.2]{W2}) says
that $(G_n,K_n)$ is a Gelfand pair if and only if
the representation of $(K_n)_{_\C}$, on polynomials on $\C^n$, is
multiplicity free.  Those groups were classified by Ka\v c 
\cite[Theorem 3]{K} in another context.  Benson, Jenkins and Ratcliff
put it together for a classification of these ``irreducible Heisenberg''
Gelfand pairs $(G_n,K_n)$.  See \cite[Theorem 4.6]{BJR}.  In 
Section \ref{sec9} below we'll look at some cases where 
$K_n$ need not be irreducible on $\C^n$.
\medskip

Ka\v c' list (as formulated in \cite[(13.2.5)]{W2}) is

{\small
\begin{equation} \label{kac-table}
\begin{tabular}{|r|l|l|l|l|}\hline
\multicolumn{5}{| c |}{Irreducible connected groups $K_n \subset U(n)$ 
multiplicity free on polynomials on $\C^n$} \\
\hline \hline
 & Group $K_n$ & Group $(K_n)_{_\C}$ & Acting on & Conditions on $n$ \\ \hline
1 & $SU(n)$ & $SL(n;\C)$ & $\C^n$ & $n \geqq 2$ \\ \hline
2 & $U(n)$ & $GL(n;\C)$ & $\C^n$ & $n \geqq 1$ \\ \hline
3 & $Sp(m)$ & $Sp(m;\C)$ & $\C^n$ & $n = 2m$ \\ \hline
4 & $U(1) \times Sp(m)$ & $\C^* \times Sp(m;\C)$ & $\C^n$ & $n = 2m$ \\ \hline
5 & $U(1) \times SO(n)$ & $\C^* \times SO(n;\C)$ & $\C^n$ & $n \geqq 2$
        \\ \hline
6 & $U(m)$ & $GL(m;\C)$ & $S^2(\C^m)$ & $m \geqq 2, \
        n = \tfrac{1}{2}m(m+1)$  \\ \hline
7 & $SU(m)$ & $SL(m;\C)$ & $\Lambda^2(\C^m)$ & $m$ odd,
        $n = \tfrac{1}{2}m(m-1)$  \\ \hline
8 & $U(m)$ & $GL(m;\C)$ & $\Lambda^2(\C^m)$ &
        $n = \tfrac{1}{2}m(m-1)$  \\ \hline
9 & $SU(\ell) \times SU(m)$ & $SL(\ell;\C) \times SL(m;\C)$ &
        $\C^\ell \otimes \C^m$ & $n = \ell m, \ \ell \ne m$ \\ \hline
10 & $U(\ell) \times SU(m)$ & $GL(\ell;\C) \times SL(m;\C)$ &
        $\C^\ell \otimes \C^m$ & $n = \ell m$ \\ \hline
11 & $U(2) \times Sp(m)$ & $GL(2;\C) \times Sp(m;\C)$ &
        $\C^2 \otimes \C^{2m}$ & $n = 4m$ \\ \hline
12 & $SU(3) \times Sp(m)$ & $SL(3;\C) \times Sp(m;\C)$ &
        $\C^3 \otimes \C^{2m}$ & $n = 6m$ \\ \hline
13 & $U(3) \times Sp(m)$ & $GL(3;\C) \times Sp(m;\C)$ &
        $\C^3 \otimes \C^{2m}$ & $n = 6m$ \\ \hline
14 & $U(4) \times Sp(4)$ & $GL(4;\C) \times Sp(4;\C)$ &
        $\C^4 \otimes \C^8$ & $n = 32$ \\ \hline
15 & $SU(m) \times Sp(4)$ & $SL(m;\C) \times Sp(4;\C)$ &
        $\C^m \otimes \C^8$ & $n = 8m, \ m \geqq 3$ \\ \hline
16 & $U(m) \times Sp(4)$ & $GL(m;\C) \times Sp(4;\C)$ &
        $\C^m \otimes \C^8$ & $n = 8m, \ m \geqq 3$ \\ \hline
17 & $U(1) \times Spin(7)$ & $\C^* \times Spin(7;\C)$ & $\C^8$ &
        $n = 8$ \\ \hline
18 & $U(1) \times Spin(9)$ & $\C^* \times Spin(9;\C)$ & $\C^{16}$ &
        $n = 16$ \\ \hline
19 & $Spin(10)$ & $Spin(10;\C)$ & $\C^{16}$ & $n = 16$ \\ \hline
20 & $U(1) \times Spin(10)$ & $\C^* \times Spin(10;\C)$ &
        $\C^{16}$ & $n = 16$ \\ \hline
21 & $U(1) \times G_2$ & $\C^* \times G_{2,\C}$ & $\C^7$ & $n=7$ \\ \hline
22 & $U(1) \times E_6$ & $\C^* \times E_{6,\C}$ & $\C^{27}$ & $n=27$ \\ \hline
\end{tabular}
\end{equation}
}
Now we have the direct systems 
{\small
\begin{equation} \label{jaw-table}
\begin{tabular}{|c|l|l|l|}\hline
\multicolumn{4}{| c |}{Direct systems $\{(H_n\rtimes K_n,K_n)\}$ of Gelfand pairs,}
\\
\multicolumn{4}{| c |}{$K_n$ connected and irreducible on $\C^n$} \\
\hline \hline
 & Group $K_n$ & Acting on & Conditions on $n$ \\ \hline
1 & $SU(n)$ & $\C^n$ & $n \geqq 2$ \\ \hline
2 & $U(n)$ & $\C^n$ & $n \geqq 1$ \\ \hline
3 & $Sp(m)$ & $\C^n$ & $n = 2m$ \\ \hline
4 & $U(1) \times Sp(m)$ & $\C^n$ & $n = 2m$ \\ \hline
5a & $U(1) \times SO(2m)$ & $\C^{2m}$ & $n= 2m \geqq 2$
        \\ \hline
5b & $U(1) \times SO(2m+1)$ & $\C^{2m+1}$ & $n = 2m+1\geqq 3$
        \\ \hline
6 & $U(m)$ & $S^2(\C^m)$ & $m \geqq 2, \
        n = \tfrac{1}{2}m(m+1)$  \\ \hline
7 & $SU(m)$ & $\Lambda^2(\C^m)$ & $m$ odd,
        $n = \tfrac{1}{2}m(m-1)$  \\ \hline
8 & $U(m)$ & $\Lambda^2(\C^m)$ &
        $n = \tfrac{1}{2}m(m-1)$  \\ \hline
9 & $SU(\ell) \times SU(m)$ &
        $\C^\ell \otimes \C^m$ & $n = \ell m, \ \ell \ne m$ \\ \hline
10 & $S(U(\ell) \times U(m))$ &
        $\C^\ell \otimes \C^m$ & $n = \ell m$ \\ \hline
11 & $U(2) \times Sp(m)$ &
        $\C^2 \otimes \C^{2m}$ & $n = 4m$ \\ \hline
12 & $SU(3) \times Sp(m)$ &
        $\C^3 \otimes \C^{2m}$ & $n = 6m$ \\ \hline
13 & $U(3) \times Sp(m)$ &
        $\C^3 \otimes \C^{2m}$ & $n = 6m$ \\ \hline
15 & $SU(m) \times Sp(4)$ &
        $\C^m \otimes \C^8$ & $n = 8m, \ m \geqq 3$ \\ \hline
16 & $U(m) \times Sp(4)$ &
        $\C^m \otimes \C^8$ & $n = 8m, \ m \geqq 3$ \\ \hline
\end{tabular}
\end{equation}
}
\hskip -.15 cm
In each case the direct system $\{K_n\}$ is both strict and parabolic.
(We separated entry 5 of Table \ref{kac-table} into entries 5a and 5b of 
Table \ref{jaw-table} in order to have the parabolic property.)
\medskip

We now suppose that $\{K_n\}$ is one of the  strict parabolic direct
system, for example one of the $16$ systems given by the rows of 
Table \ref{jaw-table}.  We retain the
notation of Section \ref{sec4} for the representations, highest weights,
unitary inclusions, etc., associated to $\{K_n\}$.
\medskip

As $U(n)$ acts on $H_n = \Im\C + \C^n$ by $k: (z,v) \mapsto (z,kv)$ it
preserves the equivalence class of each of the square integrable representations
$\pi_{n,t}$ of $H_n$.  The Mackey obstruction vanishes and $\pi_{n,t}$
extends to a unitary representation $\widetilde{\pi_{n,t}}$ of 
$H_n\rtimes U(n)$ on $\cH_{n,t}$.  See \cite[Section 4]{W0} for a geometric
proof.  We will also write $\widetilde{\pi_{n,t}}$
for its restriction, the extension of $\pi_{n,t}$ to a unitary representation
of $G_n = H_n\rtimes K_n$.  
\medskip

Denote $\pi_{n,t,\lambda} = 
\widetilde{\pi_{n,t}}\otimes \widetilde{\kappa_{n,\lambda}}$ and  
write $\cH_{n,t,\lambda}$ for its representation space 
$\cH_{n,t}\otimes \cF_{n,\lambda}$.
Fix an orthonormal basis $\{u_i\}$ of $\cF_{n,\lambda}$.  Then 
$\{w[\um] \otimes u_i\}$ is a complete orthonormal set in $\cH_{n,t,\lambda}$.  Denote the matrix coefficients by
$$
f_{\ul,\um,i,j;t}(h,k) = \langle (w[\ul] \otimes u_i), 
((\widetilde{\pi_{n,t}}\otimes 
\widetilde{\kappa_{n,\lambda}})(h,k))(w[\um]\otimes u_j)\rangle.
$$
The formal degree $\deg \pi_{n,t,\lambda} = |t|^n\deg(\kappa_{n,\lambda})$,
so the $|t|^{n/2}\deg(\kappa_{n,\lambda})^{1/2} f_{\ul,\um,i,j;t}$ form
a complete orthonormal set in the space $\cE_{n,t,\lambda}
= \cH_{n,t,\lambda} \widehat{\otimes} \cH_{n,t,\lambda}^*$ of matrix 
coefficient functions.  Given a coefficient set 
$\varphi = (\varphi_{\ul,\um,i,j}(t,\lambda))$ we have the 
functions $\Phi_{n,t,\lambda,\varphi}$ on $G_n$ given by
\begin{equation}
\Phi_{n,t,\lambda,\varphi}(hk) = 
\sum_{\ul,\um,i,j} \varphi_{\ul,\um,i,j}(t,\lambda)
|t|^{n/2}\deg(\kappa_{n,\lambda})^{1/2} f_{\ul,\um,i,j;t}(h,k), \ \ \
h \in H_n, \ k \in K_n.
\end{equation}
Here $||\Phi_{n,t,\lambda,\varphi}||^2_{\cE_{n,t,\lambda}} =
\sum_{\ul,\um,i,j} |\varphi_{\ul,\um,i,j}(t,\lambda)|^2.$  
We sum the $\Phi_{n,t,\lambda,\varphi}$ to form $L^2$ functions 
$\Psi_{n,\varphi}$ on $H_n\rtimes K_n$, given by
\begin{equation}\label{b1}
\Psi_{n,\varphi}(hk) = \sum_{\kappa_{n,\lambda} \in \widehat{K_n}} 
\deg\kappa_{r,\lambda}
\left ( \int_{-\infty}^\infty \Phi_{n,t,\lambda,\varphi}(h,k) |t|^n dt\right ).
\end{equation}
As before, for $m \geqq n$ we have a $(G_n\times G_n)$--equivariant 
isometric injection of $L^2(G_n)$ into $L^2(G_m)$ given by
\begin{equation}\label{b2}
\zeta_{m,n}(\Psi_{n,\varphi}) = 
\Psi_{m,|t|^{(n-m)/2}(\deg\kappa_{n,\lambda}/\deg\kappa_{m,\lambda})^{1/2}
	\varphi}\ \ .
\end{equation}
This gives us
\begin{theorem}\label{lim-sd}
For $n > 0$ let $K_n$ be a compact connected subgroup of $\Aut(H_n)$
such that $\{K_n\}$ is a strict parabolic direct system.
Define $G_n = H_n\rtimes K_n$.
Then there is a strict direct system $\{L^2(G_n),\zeta_{r,n}\}$ of 
$L^2$ spaces, 
whose maps $\zeta_{m,n}: L^2(G_n) \to L^2(G_m)$ are
$(G_n\times G_n)$--equivariant unitary injections.  
Let $\Pi_n$ denote the left/right regular representation of 
$G_n\times G_n$ on $L^2(G_n)$.
Note that $G: = \varinjlim (G_n) = H \rtimes K$ 
where $H = \varinjlim H_n$ and $K = \varinjlim K_n$.
Thus we have a well defined Hilbert space $L^2(G): = 
\varinjlim \{L^2(G_n), \zeta_{m,n}\}$ and a natural unitary representation 
$\Pi = \varinjlim \Pi_n$ of $G \times G$ on
$L^2(G)$.  Further, $\Pi$ is the multiplicity--free
direct integral of the irreducible representations 
$\pi_{t,\lambda} := \varinjlim \pi_{n,t,\lambda}$.
\end{theorem}

Now suppose that $\{K_n\}$ is in fact one of the $16$ systems of
Table \ref{jaw-table}.  We will use its specific properties in order to
pass from the left/right representation of $G_n\times G_n$ on $L^2(G_n)$ 
to the left regular representation of $G_n$ on $L^2(G_n/K_n)$.
\medskip

Define $G_n = H_n\rtimes K_n$.
Since $(G_n, K_n)$ is a Gelfand pair with $K_n$ irreducible
on $\C^n$, Carcano's Theorem \cite{C} says that the action of $K_n$ on
the polynomial ring $\C[C^n]$ is multiplicity free, and it picks out the
right $K_n$--invariants in $L^2(G_n)$, as follows.

\begin{lemma}\label{kntriv}
Let $\kappa_{n,\lambda} \in \widehat{K_n}$.  
Define $\widetilde{\kappa_{n,\lambda}} \in \widehat{G_n}$ by 
$\widetilde{\kappa_{n,\lambda}}(hk) = \kappa_{n,\lambda}(k)$ for
$h \in H_n$ and $k \in K_n$.  Then $\pi_{n,t,\lambda} := 
\widetilde{\pi_{n,t}}\otimes \widetilde{\kappa_{n,\lambda}}$ has a
nonzero $K_n$--fixed vector if and only if $\kappa_{n,\lambda}^*$ occurs as a
subrepresentation of $\widetilde{\pi_{n,t}}|_{K_n}$, and in that
case the space of $K_n$--fixed vectors has dimension $1$.
\end{lemma}

\noindent {\bf Proof.} This is essentially the argument in 
\cite[Section 14.5A]{W0}.  Decompose $\widetilde{\pi_{n,t}}|_{K_n} = 
\sum_{\gamma \in \widetilde{K_n}} m_\gamma\ \gamma$.  Carcano's Theorem
(\cite{C}, or see \cite[Theorem 13.2.2]{W2}) says that each $m_\gamma$ is
either $0$ or $1$.  The 
$K_n$--fixed vectors of $\widetilde{\kappa} \otimes \widetilde{\pi_{n,t}}$ all
occur in $\kappa \otimes (m_{\kappa^*}\kappa^*)$, and they form a space
of dimension $m_{\kappa^*}$.  The assertion follows. \hfill $\square$
\medskip

Since $K_n$ is compact, we can view $L^2(G_n/K_n)$ as the space
of right--$K_n$--invariant functions in $L^2(G_n)$.  With
Lemma \ref{kntriv} in mind we set
$$
\widehat{K_n}^\dagger = \{\kappa_{n,\lambda} \in \widehat{K_n} \mid
\kappa_{n,\lambda}^* \text{ occurs in the space of polynomials on } \C^n\} .
$$
Given $\kappa_{n,\lambda} \in \widehat{K_n}^\dagger$ the right $K_n$--invariant
in $\C[\C^n]\otimes \cF_{n,\lambda}^*$ is the sum over a basis of the
$\kappa_{n,\lambda}$--subspace of $\C[\C^n]$ times the dual basis of
$\cF_{n,\lambda}^*$.  Normalize it to a unit vector $u_{n,t,\lambda}$
Then the (left regular) representation of $G_n$ on $L^2(G_n/K_n)$ is 
$\sum_{\kappa_{n,\lambda} \in \widehat{K_n}^\dagger} \int_{-\infty}^\infty
\widetilde{\pi_{n,t}}\otimes \widetilde{\kappa_{n,\lambda}}$ and its 
representation space is
$
\sum_{\kappa_{n,\lambda} \in \widehat{K_n}^\dagger} \int_{-\infty}^\infty
(\cH_{n,t,\lambda}\otimes u_{n,t,\lambda}\C)\, dt.
$

\begin{proposition}\label{res-inv}
If $m \geqq n$ and $\kappa_{n,\lambda} \in \widehat{K_n}^\dagger$
then $\kappa_{m,\lambda} \in \widehat{K_m}^\dagger$, and
$\kappa_{n,\lambda}$ and $\kappa_{m,\lambda}$ have the same highest
weight $\lambda$ space.  
\end{proposition}

\noindent {\bf Proof.}
The group $K_n$ acts on $\C^n$ by some representation
$\gamma_n$, so the representation of $K_n$ on polynomials of
degree $d$ is the symmetric power $S^d(\gamma_n^*)$.  Thus we can compute
the set $X_{n,d}$ of highest weights of $K_n$ on the space $P_{n,d}$ of 
polynomials of degree $d$ on $\C^n$.  Running through the $16$ cases of 
Table \ref{jaw-table} we see that $X_{n,d} \subset X_{m,d}$.
\medskip

Now let $\lambda \in X_{n,d}$.  Let $v_{n,\lambda}$ denote a (nonzero)
highest weight $\lambda$ vector for $\gk_n$ in $P_{n,d}$, and
similarly let $v_{m,\lambda}$ denote a (nonzero) highest weight $\lambda$ 
vector for $\gk_m$.  Divide up the variables of $\C^m$ to 
$\{w_1, \dots , w_n\}$ for $\C^n$ and $\{z_{n+1}, \dots , z_m\}$ for its 
complement in $\C^m$.  Express $v_{m,\lambda} = \sum_{A,B}b_{A,B}w^Az^B$
where each term has total degree $|A|+|B| = d$.  Note that $K_n$ treats
the $z_i$ as constants.  Evaluating the $z_i$ at arbitrary constant values
$C = (c_{n+1}, \dots , c_m)$ we have a highest weight $\lambda$ vector for
$\gk_n$. By Carcano's Theorem it is a multiple of $v_{n,\lambda}$. In other 
words $v_{m,\lambda}|_{\{z = C\}} = m_{_C}v_{n,\lambda}$.  The terms
$b_{A,B}w^Az^B$ with $z$--degree $|B| > 0$ yield evaluations of 
$w$--degree $|A| < d$, and cannot contribute to any $m_{_C}v_{n,\lambda}$.
Now $b_{A,B}w^Az^B = 0$ whenever $|B| > 0$.  This shows that
$v_{m,\lambda}$ is a homogeneous polynomial of degree $d$ in the $w_j$,
as is $v_{n,\lambda}$.  We conclude that $v_{m,\lambda}$ is a nonzero
multiple of $v_{n,\lambda}$.
\hfill $\square$

\begin{corollary}\label{invar-nest}
Let $m \geqq n$.  Then every $K_n$--invariant vector in $\cH_{n,t,\lambda}$
is the image of a $K_m$--invariant vector in $\cH_{m,t,\lambda}$ under the
adjoint of the unitary map 
$\cH_{n,t,\lambda} \hookrightarrow \cH_{m,t,\lambda}$.
\end{corollary}

\noindent {\bf Proof.} Retain the notation $X_{n,d}$ for those
$\lambda$ such that $\tau_{n,\lambda}$ occurs on the space $P_{n,d}$
of polynomials of degree $d$ on $\C^n$.  If $\lambda \notin X_{n,d}$
there are no nonzero $K_n$--invariant vectors in $\cH_{n,t,\lambda}$,
so the assertion is vacuous.  Now assume $\lambda \in X_{n,d}$
and choose an orthonormal basis $\{x_1, \dots , x_{q(n)}\}$ of the 
representation space for $\tau_{n,\lambda}$ in $P_{n,d}$.  According 
to Proposition \ref{res-inv} that representation space is contained in 
the representation space for $\tau_{m,\lambda}$ in $P_{m,d}$, so the 
latter has an orthonormal basis
$\{x_1, \dots , x_{q(n)}, x_{q(n)+1}, \dots , x_{q(m)}\}$.  Let
$\{x^*_1, \dots , x^*_{q(n)}\}$ and 
$\{x^*_1, \dots , x^*_{q(n)}, x^*_{q(n)+1}, \dots , x^*_{q(m)}\}$
be the corresponding dual bases of $\cL_{n,\lambda}$ and $\cL_{m,\lambda}$.
the $K_n$--invariant vectors in $\cH_{n,t,\lambda}$ are the multiples
of $\sum_1^{q(n)} x_i\otimes x_i^*$, and the $K_m$--invariant vectors in
$\cH_{m,t,\lambda}$ are the multiples of $\sum_1^{q(m)} x_i\otimes x_i^*$.
The adjoint of unitary inclusion is orthogonal projection, which sends
$\sum_1^{q(m)} x_i\otimes x_i^*$ to $\sum_1^{q(n)} x_i\otimes x_i^*$.
\phantom{XXXXXXXX} \hfill $\square$
\medskip

Combining Theorem \ref{lim-sd}, Lemma \ref{kntriv} and Corollary 
\ref{invar-nest} we have

\begin{theorem}\label{heis-case}
Let $\{(H_n\rtimes K_n,K_n)\}$ be one of the $16$ direct systems of
{\rm Table \ref{jaw-table}}.  Denote $G_n = H_n\rtimes K_n$,
$G = \varinjlim G_n$ and $K = \varinjlim K_n$.  Then the unitary direct 
system $\{L^2(G_n),\zeta_{m,n}\}$ of {\rm Theorem \ref{lim-sd}} restricts to
a unitary direct system $\{L^2(G_n/K_n), \zeta_{m,n}\}$, the 
Hilbert space $L^2(G/K) := \varinjlim \{L^2(G_n/K_n), \zeta_{m,n}\}$ 
is the subspace of
$L^2(G) := \varinjlim L^2(G_n)$ of right--$K$--invariant 
functions, and the natural unitary representation of $G$ on 
$L^2(G/K)$ is a multiplicity free direct integral of lim--irreducible 
representations.
\end{theorem}

\begin{remark}\label{res-heis-compact}{\rm
As in Remarks \ref{res-heis} and \ref{res-compact}, if $m \geqq n$ then 
the adjoint of the direct system map 
$\zeta_{m,n}: L^2(G_n/K_n) \to L^2(G_m/K_m)$ is
orthogonal projection to the image subspace, and on each 
$(G_m \times G_m)$--irreducible direct integrand it is a scalar
multiple of restriction of functions.  The scalar is given by the
formal degree, so it depends on the integrand.} \hfill $\diamondsuit$
\end{remark}

\section{Extension to Certain Classes of Nilpotent Groups}\label{sec6}
\setcounter{equation}{0}

Theorem \ref{heis-inject} depends on four basic facts.  First, the 
$\pi_{n,t}$ are determined by their central character.  Second, 
we have good models $\cH_{n,t}$ for the representation spaces, such that 
$n$ does not appear explicitly in the formulae for the actions of the group 
elements.  Third, the injections $H_n \hookrightarrow H_m$ restrict
to isomorphisms $Z_n \cong Z_m$ of the centers.  And fourth, 
we have complete information on the Plancherel measure for the $H_n$.
In this section we consider a somewhat larger class of nilpotent
direct systems that satisfy these conditions.  
\medskip

We will need the theory of square integrable representations of nilpotent
Lie groups (\cite{MW}, or see \cite[Section 14.2]{W2} for a short exposition).
Let $N$ be a connected, simply connected nilpotent Lie group and $\gn$ its Lie
algebra.  Decompose $\gn = \gz + \gv$ and $N = Z\exp(\gv)$ where $\gz$
is the center of $\gn$.  Then $Z = \exp(\gz)$ is the center of $N$.  We
say that an irreducible unitary representation $\pi$ of $N$ is
{\sl square integrable} if its coefficient functions $f_{u,v}(g) =
\langle u, \pi(g)v\rangle$ satisfy $|f_{u,v}| \in L^2(N/Z)$.  In that case 
$\pi$ is
determined by its central character, $\pi = \pi_t$ where $t \in \gz^*$
and the central character is $\exp(\zeta) \mapsto e^{it(\zeta)}$.  In
terms of geometric quantization, $\pi_t$ corresponds to the coadjoint orbit
in $\gn^*$ consisting of all linear functionals on $\gn$ whose restriction 
to $\gz$ is $t$.  Further, the 
antisymmetric bilinear form $b_t(\xi,\eta) = t([\xi,\eta])$ on $\gv$
is nondegenerate, and (up to a positive constant that depends only on
normalizations of Haar measures) the formal degree of $\pi_t$
is $|\Pf(b_t)|$, where $\Pf(b_t)$ is the Pfaffian\footnote{Strictly speaking,
$\Pf(b_t)$ depends on a choice of basis of $\gv$, for a basis change
of determinant $a_t$ multiplies $\det b_t|_{\gv \times \gv}$ by
$\det a_t^2$ and multiplies $\Pf(b_t)$ by $\det a_t$.} of 
$b_t: \gv \times \gv \to \R$.  In fact, if $\pi_s$ is the representation
of $N$ that corresponds to $\Ad^*(N)s \subset \gn^*$, then 
$\pi_s$ is square integrable if and only if $\Pf(b_{s|_{\gz^*}}) \ne 0$.
In any case, $\Pf(b_t)$ is a polynomial function of $t$, and (again up
to a constant that depends on normalizations) $|\Pf(b_t)|$ is the
Plancherel density.  It follows that if one irreducible unitary 
representation of $N$ is square integrable then Plancherel--almost--all are.  
In the case of the Heisenberg group $H_n$, where we identified $\gz^*$ 
with $\R$, the Pfaffian corresponding to $\pi_{n,t}$ is $t^n$.  
\medskip

The point of this, from the viewpoint of commutative spaces, is that many
Gelfand pairs are of the form $(N\rtimes K,K)$ where $N$ is a connected
simply connected Lie group, $K$ is a compact subgroup of $\Aut(N)$, and
$N$ has square integrable representations.  See \cite[Theorem 14.4.3]{W2}.
In quite a few cases the groups $N$ of \cite[Theorem 14.4.3]{W2} fall 
naturally into direct systems for which we can apply the techniques of 
Section \ref{sec3}.  This is simplified by the $2$--step Nilpotent Theorem
\cite[Theorem 13.1.1]{W2} of Benson-Jenkins-Ratcliff and Vinberg, which says
that $N$ must be abelian or $2$--step nilpotent.  In a certain
sense representations treat those groups as Heisenberg groups:  

\begin{lemma} \label{is-heis1} {\rm (\cite[Lemma 14.4.1]{W2})}
Let $N$ be a connected simply connected $2$--step nilpotent
Lie group with $1$--dimensional center.  Then $N$ is isomorphic to the
Heisenberg group $H_n$ where $n = \tfrac{1}{2}(\dim_{_\R} \gn - 1)$, and in
particular $N$ has square integrable representations.
\end{lemma}

\begin{proposition} \label{is-heis2} {\rm (\cite[Proposition 14.4.2]{W2})}
Let $N$ be a connected simply connected $2$--step nilpotent
Lie group.  Let $f \in \gn^*$ such that $f|_\gz \ne 0$.
Denote $\gw_f = \{z \in \gz \mid f(z) = 0\}$ and $W_f := \exp(\gw_f)$.
Then

{\rm 1.} $W_f$ is a closed subgroup of $Z$, hence a closed normal subgroup
of $N$.

{\rm 2.} The functional $f$ is the pullback of a linear functional
$f' \in (\gn/\gw_f)^*$ and is nonzero on the central subalgebra $\gz/\gw_f$
of $\gn/\gw_f$.

{\rm 3.} The representation $[\pi_f]$ is the pullback to $N$ of the
class $[\pi_{f'}] \in \widehat{N/W_f}$.

{\rm 4.} If the representation $[\pi_f]$ is square integrable then
$[\pi_{f'}]$ is square integrable, and in that case
$N/W_f$ has center $Z/W_f$ and is isomorphic to a Heisenberg group
$H_n$ where $n = \tfrac{1}{2}\dim(\gn/\gz)$.
\end{proposition}

We now consider a strict direct system $\{N_n\}$ of $2$--step nilpotent
connected, simply connected Lie groups that have square integrable
representations, where the inclusions $\gn_n \to \gn_m$ map 
the center $\gz_n \hookrightarrow \gz_m$ and the complement 
$\gv_n \hookrightarrow \gv_m$ in decompositions $\gn_n = \gz_n + \gv_n$.
Then the direct limit algebra $\gn := \varinjlim \gn_n$ has center
$\gz := \varinjlim \gz_n$ and $\gn = \gz + \gv$ where 
$\gv = \varinjlim \gv_n$.  On the group
level, $Z = \varinjlim Z_n = \exp(\gz)$ is the center of
$N := \varinjlim N_n$ and we have $N = Z\exp(\gv)$.
\medskip

We further assume that the dimensions $\dim Z_n$ of the centers are
bounded.  Since they are non--decreasing we may assume that they are
eventually constant.  Passing to a cofinal sequence,
\begin{equation}
\gn_n \hookrightarrow \gn_m \text{ maps } \gz_n \cong \gz_m.
\end{equation}
Under that identification we write $\gz$ for all the $\gz_n$,\,\,
$\gz^*$ for all the $\gz_n^*$,\,\, and $Z$ for all the $Z_n$.
\medskip

Let $t \in \gz^*$.  Write $b_{n,t}$ for the
bilinear form  $(\xi,\eta) \mapsto t([\xi,\eta])$ on
$\gv_n$.  Then $t$ corresponds to a square integrable representation 
$\pi_{n,t}$ of $N_n$ just when the Pfaffian $\Pf(b_{n,t}) \ne 0$.
For purposes of comparing the Pfaffians as $n$ varies, we note that
$\Pf(b_{n,t})$ is specified by $t$ and a basis of $\gv_n$, so we 
simply assume that these bases are nested in the sense that the basis
of $\gv_{n+1}$ consists the basis of $\gv_n$ together with some elements
that are $b_{n+1,t}$--orthogonal to $\gv_n$.  Thus, if 
$\Pf(b_{n+1,t}) \ne 0$ then $\Pf(b_{n,t}) \ne 0$.  The converse fails 
in general, but the following lemma deals with the possibility that 
$\Pf(b_{n,t}) \ne 0 = \Pf(b_{m,t})$.  It depends on the fact \cite{MW} 
that each $\Pf(b_{n,t})$ is a polynomial function on $\gz^*$.

\begin{lemma}\label{v-basis}
Let $\ga_n \in \gz^*$ denote the zero set of $\Pf(b_{n,t})$ and 
set $\ga = \bigcup \ga_n$.  Then $\ga$ is a set of Lebesgue measure
zero in $\gz^*$.
\end{lemma}

\noindent {\bf Proof.}  Since $N_n$ has square integrable representations,
the Pfaffian $\Pf(b_{n,t})$ is a nontrivial polynomial function of
$t \in \gz^*$, so $\ga_n$ is a finite union of lower--dimensional
subvarieties of $\gz_n^*$.  Now the set $\ga$ is a countable union of 
sets $\ga_n$ of Lebesgue measure zero.  \hfill $\square$
\medskip

For convenience we define 
\begin{equation}\label{T}
T = 
\{t \in \gz^* \mid \text{ each } \Pf(b_{n,t}) \ne 0\} = \gz^* \setminus \ga.
\end{equation}
By construction, for every $t \in T$ and every index $n$ we have
a square integrable representation $\pi_{n,t} \in \widehat{N_n}$.
\medskip

Fix $t \in T$.  Then we have the hyperplane $\gw_t$ in $\gz$, and
$W_t := \exp(\gw_t)$ is a closed subgroup of $Z$.  Lemma \ref{is-heis1}
and Proposition \ref{is-heis2} tell us that each quotient $N_n/W_t$
is isomorphic to a Heisenberg group $H_{d(n)}$ and that $\pi_{n,t}$ factors
through to the square integrable representation of $N_n/W_t$ with 
central character $e^{it}$.  Now the various (as $t$ varies in $T$) 
$\pi_{n,t}$ act on the same Fock space $\cH_{d(n),t}$, 
$d(n) = \tfrac{1}{2}\dim \gv_n$, by formulae independent of $d(n)$.  
\medskip

We normalize the inner products on the $\cH_{d(n),t}$ as before,
so the $w[\um]$ form a complete orthonormal set, and realize the
space $\cE_{n,t} = \cH_{d(n),t}\widehat{\otimes} \cH_{d(n),t}^*$ of
matrix coefficients as the span of the $f_{\ul,\um,t}
: g \mapsto \langle w[\ul],\pi_{n,t}(g)w[\um]\rangle$
as in Section \ref{sec3}.  The orthogonality relations say that the
inner product on $\cE_{n,t}$ is given by 
$\langle f_{\ul,\um,t}, f_{\ul',\um',t}\rangle
= |\Pf(b_{n,t})|^{-1}$ if $\ul = \ul'$ and 
$\um = \um'$, and is $0$ otherwise.  Now the $|\Pf(b_{n,t})|^{1/2}f_{\ul,\um,t}$
form a complete orthonormal set in $\cE_{n,t}$, and as before
$\cE_{n,t}$ consists of the functions $\Phi_{n,t,\varphi}$ on $H_{d(n)}$
given by 
\begin{equation}
\Phi_{n,t,\varphi}(h) = 
\sum_{\ul,\um} \varphi_{\ul,\um}(t) |\Pf(b_{n,t})|^{1/2}f_{\ul,\um;t}(h)
\text{ with } \sum_{\ul,\um} |\varphi_{\ul,\um}(t)|^2 < \infty.  
\end{equation} 
Now
$L^2(N_n)$ is the direct integral $\int_{\gz_n^*} \cE_{n,t} |\Pf(b_{n,t})| 
dt = \int_T \cE_{n,t} |\Pf(b_{n,t})| dt$.  It consists of all functions
$\Psi_{n,\varphi}$ defined by
\begin{equation}\label{c1}
\Psi_{n,\varphi}(h) =  \int_{\gz_n^*} \Phi_{n,t,\varphi}(h) 
|\Pf(b_{n,t})| dt
= \int_T \left ( {\sum}_{\ul,\um}\varphi_{\ul,\um}(t)
	|\Pf(b_{n,t})|^{1/2} f_{\ul,\um,t} \right ) |\Pf(b_{n,t})| dt
\end{equation}
such that the functions $\varphi_{\ul,\um}: \gz_n^* \to \C$ 
are measurable with $\sum_{\ul,\um}
|\varphi_{\ul,\um}(t)|^2 < \infty$ for almost all $t \in T$ and
$\sum_{\ul,\um} |\varphi_{\ul,\um}(t)|^2 \in L^1(\gz_n^*, |\Pf(b_{n,t})| dt)$.
The norms are
\begin{equation}\label{c2}
\begin{aligned}
||\Psi_{n,\varphi}||^2_{L^2(N_n)} 
	&= \int_T ||\Phi_{n,t,\varphi}||^2_{\cE_{n,t}} 
		|\Pf(b_{n,t})| dt \\
	&= \int_T \left ( {\sum}_{\ul,\um}
		|\varphi_{\ul,\um}(t)|^2 \right ) |\Pf(b_{n,t})| dt
	= {\sum}_{\ul,\um} 
	 ||\varphi_{\ul,\um}||^2_{L^2(\gz^*,|\Pf(b_{n,t})|dt)}\ \ .
\end{aligned}
\end{equation}
As in the Heisenberg group case, the left/right representation of 
$N_n\times N_n$ on $\cE_{n,t}$ is the exterior tensor product
$\pi_{n,t} \boxtimes \pi_{n,t}^*$; it is irreducible and the left/right
regular representation of $N_n\times N_n$ on $L^2(N_n)$ is
$\Pi_n = \int_{\gz^*} (\pi_{n,t} \boxtimes \pi_{n,t}^*) |\Pf(b_{n,t})|dt$.
The argument of Lemma \ref{regrephn} goes through without change, proving

\begin{lemma}\label{regrepnn}
The left/right regular representation $\Pi_n$ of $N_n\times N_n$ on
$L^2(N_n)$ is a multiplicity free direct integral of the irreducible
unitary representations $\pi_{n,t}$.
\end{lemma}

\begin{remark} {\rm From the considerations just described, one sees that 
Lemma \ref{regrepnn} holds for every $2$--step nilpotent
Lie group that has square integrable representations.}
\end{remark}

We now continue the argument as in the Heisenberg group case.  
Suppose that the index $m \geqq n$.  Then 
$|\Pf(b_{n,t})|^{1/2} f_{\ul,\um,t} \mapsto |\Pf(b_{m,t})|^{1/2} f_{\ul,\um,t}$ 
defines an isometric injection 
$\Phi_{n,\varphi}(t) \mapsto \Phi_{m,\varphi}(t)$ of
$\cE_{n,t}$ into $\cE_{m,t}$.  The norm computation just above 
gives
\begin{equation}\label{c3}
\begin{aligned}
||\Psi_{m,|\Pf(b_{n,t})/\Pf(b_{m,t})|^{1/2}\,\varphi}||^2_{L^2(N_m)} 
&= \int_T \left ( {\sum}_{\ul,\um} |(\Pf(b_{n,t})/\Pf(b_{m,t})|\,\,
          |\varphi_{\ul,\um}(t)|^2 \right ) |\Pf(b_{m,t})|dt \\
&= \int_T \left ( {\sum}_{\ul,\um} 
	  |\varphi_{\ul,\um}(t)|^2 \right ) |\Pf(b_{n,t})|dt 
= ||\Psi_{n,\varphi}||^2_{L^2(N_n)}.
\end{aligned}
\end{equation}
Thus we have an $(N_n\times N_n)$--equivariant isometric injection 
\begin{equation}\label{c4}
\zeta'_{m,n}: L^2(N_n) \to L^2(N_m) \text{ defined by }
	\zeta'_{m,n}(\Psi_{n,\varphi}) =
	  \Psi_{m,|\Pf(b_{n,t})/\Pf(b_{m,t})|^{1/2}\varphi}.
\end{equation}
On the level of coefficients it is given by
$\zeta'_{m,n}(\Phi_{n,t,\varphi}) =
\Phi_{m,t,|\Pf(b_{n,t})/\Pf(b_{m,t})|^{1/2}\varphi}$.  In other words 
$\zeta'_{m,n}$ sends the function $\sum_{\ul,\um} \varphi_{\ul,\um}(t)
                |\Pf(b_{n,t})|^{1/2}f_{\ul,\um;t}$
on $N_n$ to the function on $N_m$ given by
$$\sum_{\ul,\um} (|\Pf(b_{n,t})/\Pf(b_{m,t})|^{1/2}
        \varphi_{\ul,\um}(t))(|\Pf(b_{m,t})|^{1/2}f_{\ul,\um;t})
= \sum_{\ul,\um} \varphi_{\ul,\um}(t) |\Pf(b_{n,t})|^{1/2}f_{\ul,\um;t}.
$$
As in the Heisenberg case we now have

\begin{theorem} \label{nilp-inj}
There is a strict direct system $\{L^2(N_n), \zeta'_{m,n}\}$ of $L^2$ spaces.
The direct system maps $\zeta'_{m,n}: L^2(N_n) \to L^2(N_m)$ are
$(N_n\times N_n)$--equivariant unitary injections.  Let $\Pi_n$ denote the
left/right regular representation of $N_n\times N_n$ on $L^2(N_n)$
and let $N = \varinjlim N_n$.  Then we have a well defined Hilbert
space $L^2(N): = \varinjlim \{L^2(N_n), \zeta'_{m,n}\}$ and a natural unitary
representation $\Pi = \varinjlim \Pi_n$ of $N \times N$ on
$L^2(N)$.  Further, that representation $\Pi$ is multiplicity--free.
\end{theorem}

\begin{remark}\label{res-nilp}{\rm
Exactly as in Remark \ref{res-heis}, the adjoint of $\zeta'_{m,n}:
L^2(N_n) \to
L^2(N_m)$ is orthogonal projection, and on each $(N_m\times N_m)$--irreducible
direct integrand of $L^2(N_m)$ it is a scalar multiple of restriction
of functions.} \hfill $\diamondsuit$
\end{remark}

\section{Structural Preliminaries}\label{sec7}
\setcounter{equation}{0}

In this section and the next we work out some structural results for a 
strict direct system $\{K_n, \varphi_{m,n}\}$ of compact connected Lie 
groups and a consistent family $\{\gamma_n\}$ of representations of the 
$K_n$ on a fixed finite dimensional vector space $\gz$.  In Section 
\ref{sec8} we will use that information to extend Theorems \ref{lim-sd} 
and \ref{heis-case} to a larger family of strict direct systems of 
nilmanifold Gelfand pairs.
\medskip

As just indicated, $\{K_n, \varphi_{m,n}\}$ is a strict direct system of
compact connected Lie groups.  Denote $K = \varinjlim \{K_n, \varphi_{m,n}\}$.
Its Lie algebra is $\gk = \varinjlim \{\gk_n, d\varphi_{m,n}\}$.
The $\gamma_n : K_n \to U(\gz)$ are unitary representations of the $K_n$ on
the finite dimensional vector space $\gz$.  They are consistent in the sense 
that $\gamma_m \cdot \varphi_{m,n} = \gamma_n$, so they define the
direct limit representation $\gamma = \varinjlim \gamma_n$ of $K$ on $\gz$.
\medskip

Let $U_n = \gamma_n(K_n)$.  The $U_n$ form an increasing sequence of compact 
connected subgroups of dimension $\leqq (\dim\gz)^2$ in the unitary 
group $U(\gz)$, so from some index on they are all the same compact connected
subgroup $U$ of $U(\gz)$.  Truncating the index set we may
assume that every $\gamma_n(K_n) = U$, in particular that, in the limit,
$U = \gamma(K)$.  Let $K_n^\dagger$ denote the 
identity component of the kernel of $\gamma_n$.  Then 
$\varphi_{m,n}(K_n^\dagger) \subset K_m^\dagger$, so we have
$K^\dagger = \varinjlim \{K_n^\dagger, \varphi_{m,n}|_{K_n^\dagger}\}$,
and $K^\dagger$ is the identity component of the kernel of $\gamma$.
In particular its Lie algebra 
$\gk^\dagger = \varinjlim \{\gk_n^\dagger, d\varphi_{m,n}|_{\gk_n^\dagger}\}$
is the kernel of $d\gamma: \gk \to \gu$.
\medskip

Since $K_n$ is compact and connected, and $K_n^\dagger$ is a closed
connected normal subgroup, $K_n$ has another closed connected normal
subgroup $L_n$ such that $K_n$ is locally isomorphic to the direct
product $K^\dagger_n \times L_n$.  On the Lie algebra level, 
$\gk_n = \gk_n^\dagger \oplus \gl_n$, direct sum of ideals.  
The semisimple part, $\gl'_n
= [\gl_n,\gl_n]$, is the direct sum of all the simple ideals of $\gk_n$
that are not contained in $\gk_n^\dagger$.  Thus $\gl'_n$ is independent
of the choice of $L_n$, and $d\varphi_{m,n}(\gl'_n) = \gl'_m$.
\medskip

\begin{proposition}\label{choose-center}
One can choose the groups $L_n$ so that $\varphi_{m,n}(L_n) = L_m$
for $m \geqq n \ggg 0$.
\end{proposition}

\noindent {\bf Proof.}  We argue by induction on $r = \dim (\gu/[\gu,\gu])$.
If $r = 0$ then the $\gl_n = \gl'_n$. so $d\varphi_{m,n}(\gl'_n) = \gl'_m$ 
says $d\varphi_{m,n}(\gl_n) = \gl_m$, and it follows that 
$\varphi_{m,n}(L_n) = L_m$.
\medskip

The group $U = \gamma(K)$ is compact and connected, and the identity component 
of its center is a torus $T$ of dimension $r$.  Suppose $r > 0$ let let
$S$ be a subtorus of dimension $r-1$ in $T$.  Now define codimension $1$ 
subgroups\, $'K_n = \gamma_n^{-1}(S) \subset K_n$ and\, 
$'K = \gamma^{-1}(S) \subset K$.  Since
$\gamma_m \cdot \varphi_{m,n} = \gamma_n$ we have $\varphi_{m,n}('K_n)
\subset {'K_m}$, and\, $'K = \varinjlim \{'K_n, \varphi_{m,n}|_{'K_n}\}$.
By induction on $r$ we have closed connected normal subgroups\, 
$'L_n \subset\, {'K_n}$ such that\, $'K_n$ isomorphic to $K_n^\dagger
\times\, 'L_n$ and $\varphi_{m,n}('L_n) =\, 'L_m$.  That gives us\,
$'L = \varinjlim\, 'L_n$ and\, $'\gl = \varinjlim\, '\gl_n$;
\, $'K$ is locally isomorphic to $K^\dagger \times\, 'L$ and
\,$'\gk = \gk^\dagger \oplus \, '\gl$.
\medskip

Let $\gw^\dagger$ denote the center of $\gk^\dagger$ and\, $'\gw$ the
center of\, $'\gl$.  Let $\gw$ denote the centralizer of\, $'\gk$ in $\gk$.  
Then $\gw$ is an abelian ideal in $\gk$ that contains $\gw^\dagger \oplus\,
'\gw$ as a subalgebra of codimension $1$.  We choose a $1$--dimensional
subalgebra\, $''\gw \subset \gw$ not contained in $\gw^\dagger \oplus\, '\gw$
and such that the $1$--parameter subgroup\, $''W = \exp(''\gw)$ 
is closed in $L$.  Define $\gl =\, {''\gw} \oplus\, {'\gl}$.  The 
corresponding analytic subgroup $L$ is a closed subgroup of $K$.  For
$n$ sufficiently large, say $n \geqq n_0$, the direct limit maps 
$\varphi_n : K_n \hookrightarrow K$ satisfy\, $''W \subset \varphi_n(K_n)$.
As indicated above, the induction hypothesis gives us\,
$'L \subset \varphi_n(K_n)$.  Thus $L \subset \varphi_n(K_n)$ for
$n \geqq n_0$.  As the $\varphi_n : K_n \hookrightarrow K$ are injective
we how have well defined closed connected subgroups $L_n = \varphi_n^{-1}(L)$
for $n \geqq n_0$, and $\gamma_n: L_n \to U$ is surjective with finite
kernel.  Thus the $\varphi_{m,n}(L_n) = L_m$, and $K_n$ is locally
isomorphic to $K_n^\dagger \times L_n$, for $m \geqq n \geqq n_0$.
That completes the proof of Proposition \ref{choose-center}.
\hfill $\square$
\medskip

For convenience of formulation we
again truncate the index set, this time 
so that $L_n = \varphi_n^{-1}(L)$ for all indices $n$.

\begin{corollary}\label{parab-cor}
Let $t \in \gz$, and let $K_{n,t}$ be its stabilizer in $K_n$.  
If one of the direct systems
$
\ \{K_n, \varphi_{m,n}\}, \ \{K_{n,t}, \varphi_{m,n}|_{K_{n,t}}\}, \
\{K_n^\dagger, \varphi_{m,n}|_{K_n^\dagger}\} \
$ is parabolic, then the other two also are parabolic.
\end{corollary}

\noindent {\bf Proof.}  Let $L_t$ denote the stabilizer of $t$ in $L$.  
Up to local isomorphism,
$$
\{K_n, \varphi_{m,n}\} = 
	\{K_n^\dagger \times L, \varphi_{m,n}|_{K_n^\dagger}\times 1\}
\text{ and }
\{K_{n,t}, \varphi_{m,n}|_{K_{n,t}}\} = 
	\{K_n^\dagger \times L_t, \varphi_{m,n}|_{K_n^\dagger}\times 1\}.
$$
In each case, the direct system is parabolic if and only if 
$\{K_n^\dagger, \varphi_{m,n}|_{K_n^\dagger}\}$ is parabolic.
\hfill $\square$
\medskip

Now Theorem \ref{lim-pw} gives us

\begin{corollary} \label{parab-cor2}
Let $t \in \gz$, and let $K_{n,t}$ be its stabilizer in $K_n$.
Suppose that the direct system $\{K_n, \varphi_{m,n}\}$ is
parabolic.  Then there are natural isometric injections
$\cF_{n,t,\lambda} \hookrightarrow \cF_{m,t,\lambda}$ for $m\geqq n$, 
from the highest weight $\lambda$ representation space of $K_{n,t}$
to that of $K_{m,t}$, and corresponding isometric injections
$\zeta''_{m,n}: f \mapsto \bigl (
(\deg \kappa_{m,t,\lambda})/(\deg\kappa_{n,t,\lambda}) \bigr )^{1/2}f$ 
on spaces of coefficient functions.
\end{corollary}

Another immediate consequence of Proposition \ref{choose-center} is

\begin{corollary}\label{compact-on-z}
Let $t \in \gz$, and let $K_{n,t}$ be its stabilizer in $K_n$.
Then $L := \varinjlim L_n$ is compact, $K = K^\dagger L$, and
$K$ is locally isomorphic to $K^\dagger \times L$.  In
particular $K$ acts on $\gz$ as a compact linear group and $\gz$ has
a $\gamma(K)$--invariant positive definite inner product.
\end{corollary}

\section{A Class of Commutative Nilmanifolds, I: Group Structure}
\label{sec8}
\setcounter{equation}{0}

In this section and the next,
we make use of the results of Sections \ref{sec6} and 
\ref{sec7} in order to extend Theorems \ref{lim-sd} and \ref{heis-case} 
to strict direct systems $\{(G_n,K_n)\}$ of Gelfand pairs that satisfy
\begin{equation}\label{nil-conditions}
\begin{aligned}
\text{(i) }&\text{each } G_n = N_n\rtimes K_n \text{ semidirect, } N_n
        \text{ connected, simply connected, }\\
        &\text{ nilpotent with square integrable representations and }
         K_n \text{ connected,}\\
\text{(ii) }&\text{the } K_n \text{ form a parabolic strict direct system},  \\
\text{(iii) }&\gn_n \hookrightarrow \gn_{n+1} \text{ maps centers }
        \gz_n \cong \gz_{n+1} \text{ and }
        \text{ complements } \gv_n \hookrightarrow \gv_{n+1}, \text{ and}\\
\text{(iv) }&\text{for each $n$ the complement $\gv_n$ is
        $\Ad(K_n)$--invariant.}
\end{aligned}
\end{equation}
We identify each $\gz_n$ with $\gz := \varinjlim \gz_n$.  Let $K_n^\dagger$ 
denote the identity component of the kernel of the action of $K_n$ on $\gz$.
The image $\Ad(K_n)|_\gz$ of that action is a compact connected group of
linear transformations of $\gz$.  Its dimension is bounded because
$\dim \gz < \infty$.  We may assume that each $\Ad(K_n)|_\gz = U$ for some 
compact connected group $U$ of linear transformations of $\gz$.
Proposition \ref{choose-center} gives us complementary closed connected normal
subgroups $L_n \subset K_n$ that map isomorphically under $K_n \hookrightarrow
K_{n+1}$, so each $L_n$ is equal to $L := \varinjlim L_n$.
Thus we have decompositions $K_n = K_n^\dagger \cdot L$ and 
$K = K^\dagger\cdot L$ where $K_n^\dagger$ is the kernel of the adjoint
action of $K_n$ on $\gz$, $K = \varinjlim K_n$, and 
$K^\dagger = \varinjlim K_n^\dagger$.  For each $n$,
$\Ad_{K_n}$ maps $L = L_n$ onto $U$ with finite kernel.  
\medskip

If $t \in \gz^*$ write $\cO_t$ for the orbit $\Ad^*(L)(t)$.  Then $\Ad^*(G)(t)
= \Ad^*(G_n)(t) = \Ad^*(K_n)(t) = \cO_t \cong L/L_t$ for each $n$.  Since
$\Ad^*(G_n)$ acts on $\cO_t$ as the compact group $L$ there is an
invariant measure $\nu_t$ derived from Haar measure on $L$; we 
normalize $\nu_t$ to total mass $1$.  Given $t \in \gz^*$ we have its
stabilizers $G_t = \{g \in G \mid \Ad^*(g)t = t\}$ and $L_t = K\cap G_t$,
and their pullbacks $G_{n,t}$ and $L_{n,t}$. 
\medskip

Let $T = \{t \in \gz^* \mid \text{ each } \Pf(b_{n,t}) \ne 0\}$, as in 
Section \ref{sec6}, and fix $t \in T$.  As in Section \ref{sec5} the
square integrable representation $\pi_{n,t}$ extends to a unitary
representation $\widetilde{\pi_{n,t}}$ of $G_{n,t} := 
N_n\rtimes K_{n,t}$ on the same representation space $\cH_{n,t}$.  If 
$\kappa_{n,t,\lambda} \in \widehat{K_{n,t}}$ has representation space 
$\cF_{n,t,\lambda}$ we write $\widetilde{\kappa_{n,t,\lambda}}$ for its
extension to a representation of $G_{n,t}$ on
$\cF_{n,t,\lambda}$ that annihilates $N_n$.  Then we have the irreducible 
unitary representation $\pi_{n,t,\lambda}^\diamondsuit := 
\widetilde{\pi_{n,t}} \otimes \widetilde{\kappa_{n,t,\lambda}}$ of
$G_{n,t}$ on $\cH_{n,t,\lambda}^\diamondsuit := 
\cH_{n,t} \otimes \cF_{n,t,\lambda}$.  
That gives us the unitary representation
$\pi_{n,t,\lambda} = \Ind_{G_{n,t}}^{G_n} (\pi_{n,t,\lambda}^\diamondsuit)$
of $G_n$.  Its representation space 
$$
\cH_{n,t,\lambda} := \int_{\cO_t} (\cH_{n,\Ad^*(k)t} 
\otimes \cF_{n,t,\lambda})\, d\nu_t(k(t))
$$ 
consists of all measurable functions
$\varphi: G_n \to \cH_{n,t,\lambda}^\diamondsuit$ such that 
$$
\begin{aligned}
\text{(i) } &\varphi(g m)=\pi_{n,t,\lambda}^\diamondsuit (m)^{-1}\varphi(\ell) 
	\text{ for } g \in G_n \text{ and }m \in G_{n,t} \text{ and }\\
\text{(ii) } &\int_{\cO_t} ||\varphi(g)||^2d\nu_t(\Ad^*(g)(t)) < \infty.
\end{aligned}
$$
In other words $\cH_{n,t,\lambda}$ is the space 
$L^2(\cO_t;\H_{n,t,\lambda}^\diamondsuit)$ of $L^2$ sections of
the homogeneous Hilbert space bundle $\H_{n,t,\lambda}^\diamondsuit \to \cO_t$
with fiber $\cH_{n,t,\lambda}^\diamondsuit$.
The action $\pi_{n,t,\lambda}$ of $G_n$ on $\cH_{n,t,\lambda}$ is
$[(\pi_{n,t,\lambda}(g))(\varphi)](g') = \varphi(g^{-1}g')$.
The inner product on $\cH_{n,t,\lambda}$ is
$\langle \varphi, \psi \rangle_{_{\cH_{n,t,\lambda}}} =
\int_{\cO_t} \langle \varphi(g), \psi(g)
\rangle_{\cH_{n,t,\lambda}^\diamondsuit}  d\nu_t(\Ad^*(g)(t))$.
\medskip

According to the Mackey little group theory, (i) $\pi_{n,t,\lambda}$ is 
irreducible, (ii) $\pi_{n,t,\lambda}$ is equivalent to
$\pi_{n,t',\lambda'}$ if and only if $t' \in \cO_t$, say
$t' = \Ad^*(\ell)t$ where $\ell \in L$, and $\Ad^*(\ell)$ carries $\lambda$
to $\lambda'$, and (iii)
Plancherel--almost--all irreducible unitary representations of
$G_n$ are of the form $\pi_{n,t,\lambda}$ where $t \in T$
and $\kappa_{n,t,\lambda} \in \widehat{K_{n,t}}$.
\medskip

Corollary \ref{parab-cor} tells us that the system $\{K_{n,t}\}$ is
parabolic.  As noted in Corollary \ref{parab-cor2} that gives us
isometric injections $\cF_{n,t,\lambda} \hookrightarrow \cF_{m,t,\lambda}$ 
of representation spaces and corresponding isometric injections of the
spaces of coefficient functions.  Those injections combine with the 
corresponding maps of Section \ref{sec6} to give us isometric injections
$\zeta''_{m,n}: f \mapsto \bigr ( (|\Pf(b_{m,t})|\deg \kappa_{m,t,\lambda})/
(|\Pf(b_{n,t})|\deg \kappa_{n,t,\lambda})\bigl )^{1/2}f$ from the space of 
coefficient functions of $\pi_{n,t,\lambda}^\diamondsuit$ to that of
$\pi_{m,t,\lambda}^\diamondsuit$.  Those come out of unitary injections 
$\cH_{n,t,\lambda}^\diamondsuit \hookrightarrow 
\cH_{m,t,\lambda}^\diamondsuit$ of the
representation spaces.  The representation space injections define 
unitary (on each fiber) injections $\H_{n,t,\lambda}^\diamondsuit
\hookrightarrow \H_{m,t,\lambda}^\diamondsuit$
of the corresponding homogeneous Hilbert space bundles over the orbit $\cO_t$.
Spaces of $L^2$ sections correspond by 
$$
L^2(\cO_t; \H_{n,t,\lambda}^\diamondsuit)
= \{\varphi \in L^2(\cO_t; \H_{m,t,\lambda}^\diamondsuit) \mid
\varphi(\ell(t)) \in \ell (\cH_{n,t,\lambda}^\diamondsuit)
\text{ for all } \ell \in L\}.
$$
The inner products on $\cH_{n,t,\lambda} = 
L^2(\cO_t; \H_{n,t,\lambda}^\diamondsuit)$ and on
$\cH_{m,t,\lambda} = 
L^2(\cO_t; \H_{m,t,\lambda}^\diamondsuit)$ are $G_n$--invariant, and
$G_n$ is irreducible on $L^2(\cO_t; \H_{n,t,\lambda}^\diamondsuit)$,
so there is a real scalar $c_{m,n} > 0$ such that 
$\varphi \mapsto c_{m,n}\varphi$
gives a $G_n$--equivariant isometric injection 
of $\cH_{n,t,\lambda}$ into $\cH_{m,t,\lambda}$.  Summarizing to this point,

\begin{proposition}\label{hnt1}
As just described we have a strict direct system
$\{\cH_{n,t,\lambda}\}$ based on $G_n$--equivariant isometric 
injections $\cH_{n,t,\lambda} \hookrightarrow \cH_{m,t,\lambda}$.
\end{proposition}

Now consider the spaces $\cE_{n,t,\lambda} := \cH_{n,t,\lambda} \boxtimes
\cH_{n,t,\lambda}^*$, Hilbert space completion of the space of
coefficient functions $f_{\varphi,\psi}: g \mapsto \langle \varphi ,
\pi_{n,t,\lambda}(g)\psi \rangle_{_{\cH_{n,t,\lambda}}}$ for 
$\varphi,\psi \in \cH_{n,t,\lambda}$.  The $(G_n\times G_n)$--invariant inner
product on $\cE_{n,t,\lambda}$ is
\begin{equation}\label{formaldeg}
\langle \varphi\boxtimes\psi , \varphi'\boxtimes\psi' 
\rangle_{_{\cE_{n,t,\lambda}}} = \tfrac{1}{d_{n,t,\lambda}} 
\langle \varphi,\varphi'\rangle_{_{\cH_{n,t,\lambda}}}
\overline{\langle\psi,\psi'\rangle_{_{\cH_{n,t,\lambda}}}}
\end{equation}
for some number $d_{n,t,\lambda} > 0$, which we
interpret as the formal degree $\deg\pi_{n,t,\lambda}$.  See Appendix A,
specifically Theorem A.1 below, for a discussion of this.
In any case, the
right/left action of $G_n\times G_n$ on $\cE_{n,t,\lambda}$ is an irreducible
unitary representation, and the isometric embeddings 
$\cH_{n,t,\lambda} \hookrightarrow \cH_{m,t,\lambda}$ define 
isometric embeddings $\zeta_{m,n}: f \mapsto 
(\deg\pi_{n+1,t,\lambda}/\deg\pi_{n,t,\lambda})^{1/2}f$ of
$\cE_{n,t,\lambda}$ into $\cE_{m,t,\lambda}$.  That gives us

\begin{proposition}\label{hnt2}
As just described we have a strict direct system
$\{\cE_{n,t,\lambda},\zeta_{m,n}\}$, whose spaces are the Hilbert spaces of 
coefficients of the representations $\pi_{n,t,\lambda}$, and whose maps are
$(G_n\times G_n)$--equivariant isometric embeddings.
\end{proposition}

The group $K$ acts on $\gz^*$ through its compact subgroup $L$, so $\gz^*$
has an $\Ad^*(K)$--invariant inner product.  Let $S$ be the unit sphere.  
Since $L$ and $S$ are compact there are only finitely many orbit types
$L_{s_1}, \dots , L_{s_\ell}$ of $L$ on $S$ (\cite{Yan}, or see
\cite[Theorem 1.7.25]{P}).  In other words every isotropy
subgroup of $L$ on $S$ is conjugate to exactly one of the $L_{s_i}$.
Thus every isotropy subgroup of $K$ on $S$ is conjugate to exactly
one of the $K_{s_i} = K^\dagger L_{s_i}$.
If $t \in \gz^*$ and $r \ne 0$ then the isotropy groups $K_t = K_{rt}$.
Now every isotropy subgroup of $K$ on $\gz^*\setminus \{0\}$ is conjugate
to exactly one of the $K_{s_i}$.  
\medskip

Decompose $S\cap T = S_1 \cup \dots \cup S_m$ where $S_i$ is the union of 
the orbits $\Ad^*(K)(s)$ in $S\cap T$ with isotropy $K_s = K_{s_i}$.  
Locally $S_i$ contains a smooth section to the action of $K$.  Thus there is 
a measurable section $\sigma_i : K\backslash S_i \to S_i$ to the action of $K$,
such that each isotropy $K_{n,\sigma_i(x)} = K_{n,s_i}$.  Let
$\Sigma_i = \sigma_i(K\backslash S_i)$.  According to the Mackey little group
theory, Plancherel almost all irreducible representations of $G_n$ are of
the form $\pi_{n,t,\lambda}$ where $t = rs_i$ with $1\leqq i \leqq m$ and
$r > 0$ and with $\kappa_{n,t,\lambda} \in \widehat{K_{n,t}}$.  That gives
\begin{equation}\label{l2}
L^2(G_n) = \sum_{i=1}^m \,\,  \sum_{\widehat{K_{n,s_i}}}\,\,
\int_{r=0}^\infty \, \int_{s \in \Sigma_i} \cE_{n,rs,\lambda}\,\, ds\, dr
\end{equation}
As $n$ increases we have the isometric equivariant injections
$\zeta_{m.n}: \cE_{n,t,\lambda} \to \cE_{m,t,\lambda}$ of Proposition
\ref{hnt2}.  When we form the discrete and continuous sums of (\ref{l2}), the
$\zeta_{m.n}$ act on the summands, where they fit together to define
isometric equivariant injections (which we also denote $\zeta_{m.n}$) from
$L^2(G_n)$ to $L^2(G_m)$, $m \geqq n$.  That yields the first assertion of

\begin{theorem}\label{hnt3}
Let $\{(G_n,K_n)\}$ be a strict direct system of commutative pairs that
satisfy {\rm (\ref{nil-conditions})}.  Then the isometric equivariant 
injections $\cE_{n,t,\lambda} \to \cE_{m,t,\lambda}$ of 
{\rm Proposition \ref{hnt2}} define $(G_n \times G_n)$--equivariant isometric
injections $\zeta_{m.n}: L^2(G_n) \to L^2(G_m)$.  That gives a direct 
system $\{L^2(G_n),\zeta_{m,n}\}$ of Hilbert spaces and equivariant 
isometric injections.  
Let $\Pi_n$ denote the left/right regular representation of $G_n\times G_n$ 
on $L^2(G_n)$ and let $G = \varinjlim G_n$.  Then we have a well defined 
Hilbert space $L^2(G): = \varinjlim \{L^2(G_n), \zeta_{m,n}\}$ and a natural 
unitary representation $\Pi = \varinjlim \Pi_n$ of $G \times G$ on
$L^2(G)$.  Further, that representation $\Pi$ is multiplicity--free.
\end{theorem}

\section{A Class of Commutative Nilmanifolds, II: Manifold Structure}
\label{sec9}
\setcounter{equation}{0}

We now pass from $L^2(G)$ to $L^2(G/K)$ for strict direct systems
of commutative spaces that satisfy (\ref{nil-conditions}).  Retain the
notation of Section \ref{sec8}.  The first step is
\begin{theorem} \label{iso-gelfand}
Let $t \in T$.  Then $(N_n\rtimes K_{n,t},K_{n,t})$ is a Gelfand pair.
In particular $K_{n,t}$ is multiplicity free on $\C[\gv_n]$.
\end{theorem}

\noindent
Note the similarity between the statement of Theorem \ref{iso-gelfand} and 
Yakimova's commutativity criterion 
(\cite[Theorem 1]{Y3}, or see \cite[Theorem 15.1.1]{W2}).
\medskip

\noindent {\bf Proof.}  We may assume that $t = s_i$, representing one of 
the orbit types of $L$ on $S \cap T$.
Suppose that $(N_n\rtimes K_{n,t},K_{n,t})$ is
not a Gelfand pair.  Then the commuting algebra $\cA$ for the representation 
$\lambda_t$ of $N_n\rtimes K_{n,t}$ on 
$L^2((N_n\rtimes K_{n,t})/K_{n,t})$ is not commutative.  Let 
$A_1, A_2 \in \cA$ with $A_1A_2 \ne A_2A_1$.  Note that
$L^2((N_n\rtimes K_n)/K_n)$ is the sum (over $j$) of the representation 
spaces of the $\widetilde{\lambda_j} :=
\Ind_{(N_n\rtimes K_{n,s_j})}^{(N_n\rtimes K_n)}(\lambda_{s_j})$.  In other
words it is the sum of the spaces $L^2_j((N_n\rtimes K_n)/K_n)$ where
$L^2_j((N_n\rtimes K_n)/K_n)$ consists of the 
$L^2(K_n/K_{n,t})$ functions $\varphi : N_n\rtimes K_n \to
L^2((N_n\rtimes K_{n,s_j})/K_{n,s_j})$ such that $\varphi(gh) = 
\lambda_{s_j}(h)^{-1}\varphi(g)$ for $g \in N_n\rtimes K_n$ and
$h \in N_n\rtimes K_{n,s_j}$.  Now define $\widetilde{A_1}$ and
$\widetilde{A_2}$ by
$(\widetilde{A_u}\varphi)(g) = A_u(\varphi(g))$.  Then
$\widetilde{A_u}(\varphi)(gh) = A_u(\varphi(gh))
= A_u(\lambda_{s_j}(h)^{-1}\varphi(g)) = \lambda_{s_j}(h)^{-1}(A_u(\varphi(g)))
= \lambda_{s_j}(h)^{-1}(\widetilde{A_u}\varphi)(g)$, so $\widetilde{A_u}$ is
a well defined linear transformation of $L^2_j((N_n\rtimes K_n)/K_n)$.  Further,
$[\widetilde{\lambda_{s_j}}(g)\cdot \widetilde{A_u}(\varphi)](g_1)
= (\widetilde{A_u}(\varphi))(g^{-1}g_1) = A_u(\varphi(g^{-1}g_1))
= A_u([\widetilde{\lambda_{s_j}}(g)\varphi](g_1))
= [\widetilde{A_u}(\widetilde{\lambda_{s_j}}(g)\varphi)](g_1)$, so 
$\widetilde{A_u}$ is an intertwining operator for $\widetilde{\lambda_{s_j}}$.
As the $A_u$ do not commute, neither do the $\widetilde{A_u}$.  Since
$(N_n\rtimes K_n,K_n)$ is a Gelfand pair this is a contradiction.  We
conclude that $(N_n\rtimes K_{n,t},K_{n,t})$ is a Gelfand pair.
In particular, now, $K_{n,t}$ is multiplicity free on $\C[\gv_n]$ by
Carcano's Theorem.
\hfill $\square$
\medskip

Recall the Hilbert bundle model for the induced representation
$\pi_{n,t,\lambda} \in \widehat{G_{n,t}}$ given by
$\pi_{n,t,\lambda} = \Ind_{G_{n,t}}^{G_n}(\pi_{n,t,\lambda}^\diamondsuit)$. 
The representation space $\cH_{n,t,\lambda}$ 
of $\pi_{n,t,\lambda}$ consists of all
$L^2(K_n/K_{n,t})$ sections of the homogeneous bundle 
$p: \H_{n,t,\lambda}^\diamondsuit \to G_n/G_{n,t} = K_n/K_{n,t}$ whose typical
fiber is the representation space $\cH_{n,t,\lambda}^\diamondsuit$ of
$\pi_{n,t,\lambda}^\diamondsuit$.  Given $k \in K_n$ we write
$k\cdot \cH_{n,t,\lambda}^\diamondsuit$ for the fiber $p^{-1}(kK_{n,t})$.
Let $u \in \cH_{n,t,\lambda}^\diamondsuit$ be a
$\pi_{n,t,\lambda}^\diamondsuit (K_{n,t})$--fixed unit vector.
Then $u$ belongs to the fiber $1\cdot \cH_{n,t,\lambda}^\diamondsuit$, and
$k\cdot u \in k\cdot \cH_{n,t,\lambda}^\diamondsuit$ depends only on the
coset $kK_{n,t}$.  Define a section
\begin{equation}\label{induced-invariant}
\sigma_u : K_n/K_{n,t} \to  \H_{n,t,\lambda}^\diamondsuit \text{ by }
	\sigma_u(kK_{n,t}) = k\cdot u.
\end{equation}
Then $\sigma_u$ is a $\pi_{n,t,\lambda}(K_n)$--invariant unit vector in
the Hilbert space $\cH_{n,t,\lambda}$.  (We will also
write $\varphi_u$ for the corresponding function 
$G_n \to \cH_{n,t,\lambda}^\diamondsuit$ such that
$\varphi_u(gg_t) = \pi_{n,t,\lambda}^\diamondsuit(g_t)^{-1}(\varphi_u(g))$
for $g \in G_n$ and $g_t \in G_{n,t}$.)  
Conversely if $\sigma$ is a $\pi_{n,t,\lambda}(K_n)$--invariant unit vector in
$\cH_{n,t,\lambda}$, then $\sigma(1K_{n,t}) = cu$ where $|c| = 1$ by
$K_{n,t}$--invariance, and then $\sigma = c\sigma_u$ by $K$--invariance.
In summary,
\begin{lemma}\label{k-invariant}
Let $t \in T$ and let $u$ be the unique {\rm (up to scalar multiple)} 
$\pi_{n,t,\lambda}^\diamondsuit (K_{n,t})$--fixed unit vector in 
$\cH_{n,t,\lambda}^\diamondsuit$  Then the section $\sigma_u$, given by 
{\rm (\ref{induced-invariant})}, is the unique {\rm (up to scalar multiple)}
$\pi_{n,t,\lambda}(K)$--fixed unit vector in $\cH_{n,t,\lambda}$.
\end{lemma}

By Theorem \ref{iso-gelfand} we can apply Proposition \ref{hnt2}
to the function spaces $\cE_{n,t,\lambda}^\diamond = 
\cH_{n,t,\lambda}^\diamond \boxtimes (\cH_{n,t,\lambda}^\diamond)^*$ 
on the groups $G_{n,t} = N_n\rtimes K_{n,t}$.  Now
combining Proposition \ref{hnt2} and Lemma \ref{k-invariant} we have
\begin{proposition}\label{restrict-invariants}
If orthogonal projection 
$\cE_{n+1,t,\lambda}^\diamondsuit \to \cE_{n,t,\lambda}^\diamondsuit$ sends 
a nonzero right $K_{n+1,t}$--invariant function to a nonzero
right $K_{n,t}$--invariant function, then orthogonal projection
$\cE_{n+1,t,\lambda} \to \cE_{n,t,\lambda}$ sends a nonzero right
$K_{n+1}$--invariant function to a nonzero right $K_n$--invariant function.
\end{proposition}

Vinberg (\cite{V1}, \cite{V2}; or see \cite[Table 13.4.1]{W2}) classified
the maximal irreducible nilpotent Gelfand pairs.  A Gelfand pair
$(G_n,K_n)$ is called {\sl maximal} if it is not obtained from another
Gelfand pair $(G'_n,K'_n)$ by the construction
$(G_n,K_n) = (G'_n/C,K'_n/(K'_n \cap C))$ for any nontrivial closed 
connected central subgroup $C$ of $G'_n$.  And $(G_n,K_n)$ is called 
{\sl irreducible} if $\Ad(K_n)$ is irreducible on $\gv_n = \gn_n / \gz$.
Here is Vinberg's classification of maximal irreducible nilpotent Gelfand 
pairs; see \cite{W2} for the notation.

{\normalsize
\begin{equation} \label{vin-table}
\begin{tabular}{|r|c|c|c|c|c|}\hline
\multicolumn{6}{| c |}{Maximal Irreducible Nilpotent Gelfand
        Pairs $(N_n\rtimes K_n,K_n)$ \quad (\cite{V1}, \cite{V2})}\\
\hline \hline
 & Group $K_n$ & $\gv_n$ & $\gz$ & 
   $\begin{smallmatrix} U(1) \text{ is} \\ \text{needed if}\end{smallmatrix}$ &
   $\begin{smallmatrix} \text{ max }\\\text{requires}\end{smallmatrix}$ 
   \\ \hline
1 & $SO(n)$ & $\R^n$ & $\Skew\R^{n\times n} = \gs\go(n)$ &  & \\ \hline
2 & $Spin(7)$ & $\R^8 = \O$ & $\R^7 = \Im\O$  &  & \\ \hline
3 & $G_2$ & $\R^7 = \Im\O$ & $\R^7 = \Im\O$ &  & \\ \hline
4 & $U(1)\cdot SO(n)$ & $\C^n$ & $\Im\C$ & & $n\ne 4$ \\ \hline
5 & $(U(1)\cdot) SU(n)$ & $\C^n$ & $\Lambda^2\C^n \oplus\Im\C$ &
        $n$ odd &  \\ \hline
6 & $SU(n), n$ odd & $\C^n$ & $\Lambda^2\C^n$ &  & \\ \hline
7 & $SU(n), n$ odd & $\C^n$ & $\Im\C$ &  & \\ \hline
8 & $U(n)$ & $\C^n$ & $\Im \C^{n\times n} = \gu(n)$ &  & \\ \hline
9 & $(U(1)\cdot) Sp(n)$ & $\H^n$ & $\Re \H^{n \times n}_0 \oplus \Im\H$ &  &
        \\ \hline
10 & $U(n)$ & $S^2\C^n$ & $\R$ & & \\ \hline
11 & $(U(1)\cdot) SU(n), n \geqq 3$ & ${\Lambda}^2\C^n$ & $\R$ & $n$ even &
        \\ \hline
12 & $U(1)\cdot Spin(7)$ & $\C^8$ & $\R^7 \oplus \R$ & & \\ \hline
13 & $U(1)\cdot Spin(9)$ & $\C^{16}$ & $\R$ & & \\ \hline
14 & $(U(1)\cdot) Spin(10)$ & $\C^{16}$ & $\R$ & & \\ \hline
15 & $U(1)\cdot G_2$ & $\C^7$ & $\R$ & & \\ \hline
16 & $U(1)\cdot E_6$ & $\C^{27}$ & $\R$ & & \\ \hline
17 & $Sp(1)\times Sp(n)$ & $\H^n$ & $\Im \H = \gs\gp(1)$ & & $n \geqq 2$
        \\ \hline
18 & $Sp(2)\times Sp(n)$ & $\H^{2\times n}$ &
        $\Im \H^{2\times 2} = \gs\gp(2)$ & & \\ \hline
19 & $(U(1)\cdot) SU(m) \times SU(n)$ &  &  &  & \\
   & $m,n \geqq 3$ & $\C^m\otimes \C^n$ & $\R$ & $m=n$ &   \\ \hline
20 & $(U(1)\cdot) SU(2) \times SU(n)$ & $\C^2 \otimes \C^n$ &
        $\Im \C^{2\times 2} = \gu(2)$ & $n=2$ & \\ \hline
21 & $(U(1)\cdot) Sp(2) \times SU(n)$ & $\H^2\otimes \C^n$ & $\R$ &
        $n \leqq 4$ & $n \geqq 3$ \\ \hline
22 & $U(2)\times Sp(n)$ & $\C^2 \otimes \H^n$ & $\Im \C^{2\times 2} = \gu(2)$ &
        & \\ \hline
23 & $U(3)\times Sp(n)$ & $\C^3 \otimes \H^n$ & $\R$ & & $n \geqq 2$
        \\ \hline
\end{tabular}
\end{equation}
}

As noted in \cite{W2}, in Table \ref{vin-table} one 
often can replace $K_n$ by a smaller 
group in such a way that $(G_n,K_n)$
continues to be a Gelfand pair.  For example, in Table \ref{vin-table},
Item 2, where $N_n$ is the octonionic Heisenberg group $H_{\O,1}$, the
pairs $(N_n\rtimes Spin(7),Spin(7))$,
$(N_n\rtimes Spin(6),Spin(6))$ and $(N_n \rtimes Spin(5),Spin(5))$ all are
Gelfand pairs; see \cite[Proposition 5.6]{L}.  This is the tip of the
iceberg for the classification of commutative nilmanifolds.  A systematic
analysis is given in \cite{Y2}; or see \cite[Chapter 15]{W2}.
\medskip

The strict direct systems in Table \ref{vin-table}, with $\dim \gz_n$ 
bounded, are
as follows.  Here we split entry line 4 of Table \ref{vin-table} so that
$\{K_n\}$ is parabolic, and we split entry 20 into two essentially different
cases.

{\footnotesize
\begin{equation} \label{ind-vin-table}
\begin{tabular}{|r|c|c|c|c|c|}\hline
\multicolumn{6}{| c |}{Direct Systems of Maximal Irreducible Nilpotent Gelfand
        Pairs $(N_n\rtimes K_n,K_n)$}\\
\hline \hline
 & Group $K_n$ & $\gv_n$ & $\gz_n$ & 
   $\begin{smallmatrix} U(1) \text{ is} \\ \text{needed if}\end{smallmatrix}$ &    $\begin{smallmatrix} \text{ max }\\\text{requires}\end{smallmatrix}$
   \\ \hline
4a & $U(1)\cdot SO(2n)$ & $\C^{2n}$ & $\Im\C$ & & $n\ne 2$ \\ \hline
4b & $U(1)\cdot SO(2n+1)$ & $\C^{2n+1}$ & $\Im\C$ & &  \\ \hline
7 & $SU(n), n$ odd & $\C^n$ & $\Im\C$ &  & \\ \hline
10 & $U(n)$ & $S^2\C^n$ & $\R$ & & \\ \hline
11 & $(U(1)\cdot) SU(n), n \geqq 3$ & ${\Lambda}^2\C^n$ & $\R$ & $n$ even &
        \\ \hline
17 & $Sp(1)\times Sp(n)$ & $\H^n$ & $\Im \H = \gs\gp(1)$ & & $n \geqq 2$
        \\ \hline
18 & $Sp(2)\times Sp(n)$ & $\H^{2\times n}$ &
        $\Im \H^{2\times 2} = \gs\gp(2)$ & & \\ \hline
19 & $(U(1)\cdot) SU(m) \times SU(n)$ &  &  &  & \\
   & $m,n \geqq 3$ & $\C^m\otimes \C^n$ & $\R$ & $m=n$ &   \\ \hline
20a & $SU(2) \times SU(n), n \geqq 3$ & $\C^2 \otimes \C^n$ &
        $\Im \C^{2\times 2} = \gu(2)$ &  & \\ \hline
20b & $U(2) \times SU(n)$ & $\C^2 \otimes \C^n$ &
        $\Im \C^{2\times 2} = \gu(2)$ &  & \\ \hline
21 & $(U(1)\cdot) Sp(2) \times SU(n)$ & $\H^2\otimes \C^n$ & $\R$ &
        $n \leqq 4$ & $n \geqq 3$ \\ \hline
22 & $U(2)\times Sp(n)$ & $\C^2 \otimes \H^n$ & $\Im \C^{2\times 2} = \gu(2)$ &
        & \\ \hline
23 & $U(3)\times Sp(n)$ & $\C^3 \otimes \H^n$ & $\R$ & & $n \geqq 2$
        \\ \hline
\end{tabular}
\end{equation}}
In each case of Table \ref{ind-vin-table}, \cite[Theorem 14.4.3]{W2} says that 
$N_n$ has square integrable representations.
In the cases $\dim \gz > 1$ of Table \ref{ind-vin-table} we have
$K_n = K'\cdot K_n''$ where the big factor $K_n''$ acts trivially
on $\gz$ and the small factor $K'$ acts on $\gz$ by its adjoint 
representation.  Summarizing these observations,

\begin{proposition} Each of the thirteen direct systems $\{(G_n,K_n)\}$ of
{\rm Table \ref{ind-vin-table}} has the properties {\rm (i)} $\{K_n\}$ is 
parabolic {\rm (ii)} the $\{K_{n,s_i}\}$ are parabolic and {\rm (iii)}
$N_n$ has square integrable representations.
\end{proposition}

The following result is immediate from the multiplicity free part of 
Theorem \ref{iso-gelfand} and the
argument of Lemma \ref{kntriv}.  

\begin{corollary}\label{kntriv-vin}
Let $\{(G_n,K_n)\}$ be one of the thirteen direct systems
{\rm Table \ref{ind-vin-table}} and let $t \in T$.
Let $\kappa_{n,\lambda} \in \widehat{K_{n,t}}$.
Define $\widetilde{\kappa_{n,t,\lambda}} \in \widehat{N_n \rtimes K_{n,t}}$ 
by $\widetilde{\kappa_{n,t,\lambda}}(h,k) = \kappa_{n,t,\lambda}(k)$.
Then $\pi_{n,t,\lambda}^\diamondsuit := 
\widetilde{\pi_{n,t}}\otimes \widetilde{\kappa_{n,t,\lambda}}$ has a
nonzero $K_{n,t}$--fixed vector if and only if $\kappa_{n,t,\lambda}^*$ 
occurs as a subrepresentation of $\widetilde{\pi_{n,t}}|_{K_{n,t}}$, and in 
that case the space of $K_{n,t}$--fixed vectors has dimension $1$.
\end{corollary}

Corollary \ref{kntriv-vin} lets us apply the argument of 
Corollary \ref{invar-nest} to the representation spaces 
$\cH_{n,t,\lambda}^\diamondsuit := \cH_{n,t} \otimes \cF_{n,t,\lambda}$
of the $\pi_{n,t,\lambda}^\diamondsuit$.  If $m \geqq n$ it shows that
every $K_{n,t}$--invariant vector in $\cH_{n,t,\lambda}^\diamondsuit$ is
the image of a $K_{m,t}$--invariant vector under the adjoint of the unitary
map $\cH_{n,t,\lambda}^\diamondsuit \to \cH_{m,t,\lambda}^\diamondsuit$.
Combining this with Proposition \ref{hnt2} we see that orthogonal projection
$\cE_{m,t,\lambda}^\diamondsuit \to \cE_{n,t,\lambda}^\diamondsuit$
sends nonzero right $K_{m,t}$--invariant functions to nonzero right
$K_{n,t}$--invariant functions.
(Here recall the space $\cE_{n,t,\lambda}^\diamondsuit := 
\cH_{n,t,\lambda}^\diamondsuit \otimes (\cH_{n,t,\lambda}^\diamondsuit)^*$ 
of functions on $G_{n,t}$.)  Now Proposition \ref{restrict-invariants} gives
us
\begin{proposition}\label{restrict-more-invariants}
Let $\{(G_n,K_n)\}$ be one of the thirteen direct systems
{\rm Table \ref{ind-vin-table}}.  Let $t \in T$ and
$m \geqq n$.  Then orthogonal projection   
$\cE_{m,t,\lambda} \to \cE_{n,t,\lambda}$ sends nonzero right $K_m$--invariant
invariant functions to nonzero right $K_n$--invariant functions.
\end{proposition}

Combining Theorem \ref{hnt3} with Corollary \ref{kntriv-vin} and
Proposition \ref{restrict-more-invariants} we arrive at

\begin{theorem}\label{big-hnt}
Let $\{(G_n,K_n)\}$ be one of the thirteen direct systems of
{\rm Table \ref{ind-vin-table}}.  Denote $G = \varinjlim G_n$ and
$K = \varinjlim K_n$.  Then the unitary direct system $\{L^2(G_n)\}$
of {\rm Theorem \ref{hnt3}} restricts to a unitary direct system
$\{L^2(G_n/K_n)\}$, the
Hilbert space $L^2(G/K) := \varinjlim L^2(G_n/K_n)$ is the subspace of
$L^2(G) := \varinjlim L^2(G_n)$ consisting of right--$K$--invariant
functions, and the natural unitary representation of $G$ on
$L^2(G/K)$ is a multiplicity free direct integral of lim--irreducible
representations.
\end{theorem}

Now we go past the cases that require irreducibility of $K_n$ on $\gv_n$.
\medskip

In the Table \ref{vin-table-ipms} below, $\gh_{n;\F}$ denotes the generalized
Heisenberg algebra $\Im\F + \F^n$ of real dimension $1+ n\dim_\R\F$ where 
$\F$ denotes the complex number field $\C$, the quaternion algebra 
$\H$, or the octonion algebra $\O$.  It is the Lie algebra of
the generalized Heisenberg group $H_{n;\F}$ given by
\begin{equation}
\begin{aligned}
H_{n,\F}: \text{ real vector } &\text{space } \Im\F + \F^n 
\text{ with group composition} \\
&(z,w)(z',w') = (z + z' + \Im h(w,w'), w + w')
\end{aligned}
\end{equation}
where $h$ is the standard positive definite hermitian form on $\F^n$.
The generalized Heisenberg groups $H_{n;\F}$ all have square integrable
representations \cite[Theorem 14.3.1]{W2}.

In Table \ref{vin-table-ipms} we have direct sum decompositions
\begin{equation}
\gn = \gn' \oplus \gz'' \text{ where } \gn' \text{ has center } \gz' = [\gn,\gn]
\text{ and } \gn \text{ has center } \gz' \oplus \gz'' = [\gn,\gn] \oplus \gz''.
\end{equation} 
and $K$--stable vector space decompositions
\begin{equation}
\gn = \gz + \gv \ \ \text{ and } \ \ \gn' = \gz' + \gv = [\gn,\gn] + \gv.
\end{equation}

Finally, in Table \ref{vin-table-ipms}, $\gs\gu(n)$ does not mean
the Lie algebra, but simply denotes its
underlying vector space, the space of $n\times n$ 
skew--hermitian complex matrices, as a module for $\Ad(U(n))$ or
$\Ad(SU(n))$.
\medskip

Here is a small reformulation of Yakimova's classification of 
indecomposable, principal, maximal and
$Sp(1)$--saturated commutative pairs $(N\rtimes K,K)$, where the action of
$K$ on $\gn/[\gn,\gn]$ is reducible. Compare \cite[Table 13.4.4]{W2}.
See \cite{Y2} or \cite{W2} for the technical definitions; for our purposes 
it suffices to note that these are the 
basic building blocks for the complete classification described in \cite{Y2}
and \cite{W2}.  We omit the case $[\gn,\gn] = 0$, where $N = \R^n$
and $K$ is any closed subgroup of the orthogonal group $O(n)$.

{\footnotesize
\begin{equation} \label{vin-table-ipms}
\begin{tabular}{|r|l|l|l|l|l|}\hline
\multicolumn{6}{| c |}{Maximal Indecomposable Principal Saturated Nilpotent
Gelfand Pairs $(N\rtimes K,K)$,}\\
\multicolumn{6}{| c |}{$N$ Nonabelian Nilpotent, Where the Action of $K$ on
$\gn/[\gn,\gn]$ is Reducible}\\
\hline \hline
 & Group $K$ & $K$--module $\gv$ & module $\gz' = [\gn,\gn]$ &
	module $\gz''$ & Algebra $\gn'$
        \\ \hline
1 & $U(n)$ & $\C^n$ & $\R$ & $\gs\gu(n)$ &
        $\gh_{n;\C}$  \\ \hline
2 & $U(4)$ & $\C^4$ & $\Im\C \oplus \Lambda^2\C^4$ & $\R^6$ &
        $\Im\C + \Lambda^2\C^4 +\C^4$ \\ \hline
3 & $U(1)\times U(n)$ & $\C^n \oplus \Lambda^2\C^n$ & $\R \oplus \R$ & $0$ &
        $\gh_{n;\C} \oplus \gh_{n(n-1)/2;\C}$ \\ \hline
4 & $SU(4)$ & $\C^4 = \H^2$ & $\Im\C \oplus \Re\H^{2\times 2}$ &
        $\R^6$ & $\Im\C + \Re\H^{2\times 2} + \C^4$ \\ \hline
5 & $U(2)\times U(4)$ & $\C^{2\times 4}$ & $\Im \C^{2\times 2}$ & $\R^6$ &
        $\Im \C^{2\times 2} + \C^{2\times 4}$ \\ \hline
6 & $S(U(4)\times U(m))$ & $\C^{4 \times m}$ & $\R$ & $\R^6$ &
        $\gh_{4m;\C}$ \\ \hline
7 & $U(m)\times U(n)$ & $\C^{m\times n}\oplus \C^m$ & $\R \oplus \R$ & $0$ &
        $\gh_{mn;\C} \oplus \gh_{m;\C}$ \\ \hline
8 & $U(1)\times Sp(n)\times U(1)$ & $\C^{2n}\oplus \C^{2n}$ & $\R\oplus \R$ &
         $0$ & $\gh_{2n;\C} \oplus \gh_{2n;\C}$ \\ \hline
9 & $Sp(1)\times Sp(n)\times U(1)$ & $\H^n \oplus \H^n$ & $\Im \H \oplus \R$ &
         $0$ & $\gh_{n;\H} \oplus \gh_{2n;\C}$ \\ \hline
10 & $Sp(1)\times Sp(n)\times Sp(1)$ & $\H^n\oplus\H^n$ & $\Im\H\oplus\Im\H$ &
         $0$ & $\gh_{n;\H} \oplus \gh_{n;\H}$ \\ \hline
11 & $Sp(n)\times\{Sp(1),U(1),\{1\}\}$ & $\H^n$ & $\Im\H$ &
        $\H^{n\times m}$ & $\gh_{n;\H}$ \\
   & \phantom{XXXXX} $\times Sp(m)$ & & & & \\ \hline
12 & $Sp(n)\times\{Sp(1),U(1),\{1\}\}$ & $\H^n$ &
	$\Re \H^{n\times n}_0$ &
        $\Im \H$ & $\gh_{n;\H}$ \\ \hline
13 & $Spin(7)\times \{SO(2),\{1\}\}$ & $\R^8=\O$ &
        $\R^7 = \Im\O$ & $\R^{7\times 2}$ &
        $\gh_{1;\O}$ \\ \hline
14 & $U(1)\times Spin(7)$ & $\C^7$ & $\R$ & $\R^8$ &
        $\gh_{7;\C}$ \\ \hline
15 & $U(1)\times Spin(7)$ & $\C^8$ & $\R$ & $\R^7$ &
        $\gh_{8;\C}$ \\ \hline
16 & $U(1)\times U(1)\times Spin(8)$ & $\C^8_+\oplus\C^8_-$ & $\R\oplus\R$ &
        $0$ & $\gh_{8;\C} \oplus \gh_{8;\C}$ \\ \hline
17 & $U(1)\times Spin(10)$ & $\C^{16}$ & $\R$ & $\R^{10}$ &
        $\gh_{16;\C}$  \\ \hline
18 & $\{SU(n),U(n),U(1)Sp(\tfrac{n}{2})\}$ &
        $\C^{n\times 2}$ & $\R$ & $\gs\gu(2)$ &
        $\gh_{2n;\C}$ \\
   & \phantom{XXXXX} $\times SU(2)$ & & & & \\ \hline
19 & $\{SU(n),U(n),U(1)Sp(\tfrac{n}{2})\}$ &
        $\C^{n\times 2}\oplus\C^2$ & $\R\oplus\R$ & $0$ &
        $\gh_{2n;\C} \oplus \gh_{2;\C}$ \\
   & \phantom{XXXXX} $\times U(2)$ & & & & \\ \hline
   & $\{SU(n),U(n),U(1)Sp(\tfrac{n}{2})\}$ & & & & \\
20 & \phantom{XXXXX}$\times SU(2)\times$ &
        $\C^{n\times 2}\oplus\C^{2\times m}$ &
        $\R \oplus \R$ & $0$ & $\gh_{2n;\C} \oplus \gh_{2m;\C}$ \\
   & $\{SU(m),U(m),U(1)Sp(\tfrac{m}{2})\}$ & & & & \\ \hline
21 & $\{SU(n),U(n),U(1)Sp(\tfrac{n}{2})\}$ &
        $\C^{n\times 2}\oplus\C^{2\times 4}$ &
        $\R \oplus \R$ & $\R^6$ & 
	$\gh_{2n;\C} \oplus \gh_{8;\C}$ \\
   & \phantom{XXXXX} $\times SU(2) \times U(4)$ & & & & \\ \hline
22 & $U(4)\times U(2)$ & $\C^{4\times 2}$ &
        $\R$ & $\R^6 \oplus \gs\gu(2)$ & 
	$\gh_{8;\C}$ \\ \hline
23 & $U(4)\times U(2)\times U(4)$ & $\C^{4\times 2} \oplus \C^{2\times 4}$ &
        $\R\oplus\R$ & $\R^6 \oplus \R^6$ & 
	$\gh_{8;\C} \oplus \gh_{8;\C}$ \\ \hline
24 & $U(1)\times U(1)\times SU(4)$ & $\C^4\oplus\C^4$ &
        $\R\oplus\R$ & $\R^6$ & 
	$\gh_{4;\C} \oplus \gh_{4;\C}$ 
	\\ \hline
25 & $(U(1)\cdot)SU(4)(\cdot SO(2))$ & $\C^4$ &
        $\R^{6\times 2}$ & 
	$\R$ & $\gh_{4;\C}$ \\ \hline
\end{tabular}
\end{equation}
}

In each case of Table \ref{vin-table-ipms}, 
the group $N = N' \times Z''$ has square integrable
representations \cite[Theorem 14.3.1]{W2}.  In fact, if $t \in \gz^*$
we decompose $t = t' + t''$ where $t'(\gz'') = 0 = t''(\gz')$, and
then $\Pf(b_t) = \Pf(b_{t'})$, independent of $t''$.
\medskip

The strict direct systems in Table \ref{vin-table-ipms}, with $\dim \gz'_n$
bounded, are as follows.  Here the index $\ell$ can be $n$ or
$(m,n)$, the group $G_\ell = N_\ell \rtimes K_\ell$, and the subgroup
$G'_\ell = N'_\ell \rtimes K_\ell$.
\medskip

{\footnotesize
\begin{equation} \label{ind-vin-table-ipms}
\begin{tabular}{|c|l|l|l|l|l|}\hline
\multicolumn{6}{| c |}{Strict Direct Systems $\{(G_\ell,K_\ell)\}$ 
and $\{(G'_\ell,K'_\ell)\}$
of Gelfand Pairs From Table \ref{vin-table-ipms} with 
$\dim \gz'_\ell$ Bounded}\\
\hline \hline
 & Group $K_\ell$ & $K_\ell$--module $\gv_\ell$ & module $\gz_\ell' = [\gn_\ell,\gn_\ell]$ &
	module $\gz''_\ell$ & Algebra $\gn'_\ell$
        \\ \hline
1 & $U(n)$ & $\C^n$ & $\R$ & $\gs\gu(n)$ &
        $\gh_{n;\C}$  \\ \hline
3 & $U(1)\times U(n)$ & $\C^n \oplus \Lambda^2\C^n$ & $\R \oplus \R$ & $0$ &
        $\gh_{n;\C} \oplus \gh_{n(n-1)/2;\C}$ \\ \hline
6 & $S(U(4)\times U(m))$ & $\C^{4 \times m}$ & $\R$ & $\R^6$ &
        $\gh_{4m;\C}$ \\ \hline
7 & $U(m)\times U(n)$ & $\C^{m\times n}\oplus \C^m$ & $\R \oplus \R$ & $0$ &
        $\gh_{mn;\C} \oplus \gh_{m;\C}$ \\ \hline
8 & $U(1)\times Sp(n)\times U(1)$ & $\C^{2n}\oplus \C^{2n}$ & $\R\oplus \R$ &
         $0$ & $\gh_{2n;\C} \oplus \gh_{2n;\C}$ \\ \hline
9 & $Sp(1)\times Sp(n)\times U(1)$ & $\H^n \oplus \H^n$ & $\Im \H \oplus \R$ &
         $0$ & $\gh_{n;\H} \oplus \gh_{2n;\C}$ \\ \hline
10 & $Sp(1)\times Sp(n)\times Sp(1)$ & $\H^n\oplus\H^n$ & $\Im\H\oplus\Im\H$ &
         $0$ & $\gh_{n;\H} \oplus \gh_{n;\H}$ \\ \hline
11a & $Sp(n)\times Sp(1)\times Sp(m)$ & $\H^n$ & $\Im\H$ &
        $\H^{n\times m}$ & $\gh_{n;\H}$ \\ \hline
11b & $Sp(n)\times U(1)\times Sp(m)$ & $\H^n$ & $\Im\H$ &
        $\H^{n\times m}$ & $\gh_{n;\H}$ \\ \hline
11c & $Sp(n)\times \{1\} \times Sp(m)$ & $\H^n$ & $\Im\H$ &
        $\H^{n\times m}$ & $\gh_{n;\H}$ \\ \hline
18a & $SU(n) \times SU(2)$ &
        $\C^{n\times 2}$ & $\R$ & $\gs\gu(2)$ &
        $\gh_{2n;\C}$ \\ \hline
18b & $U(n) \times SU(2)$ &
        $\C^{n\times 2}$ & $\R$ & $\gs\gu(2)$ &
        $\gh_{2n;\C}$ \\ \hline
18c & $U(1)Sp(\tfrac{n}{2}) \times SU(2)$ &
        $\C^{n\times 2}$ & $\R$ & $\gs\gu(2)$ &
        $\gh_{2n;\C}$ \\ \hline
19a & $SU(n) \times U(2)$ &
        $\C^{n\times 2}\oplus\C^2$ & $\R\oplus\R$ & $0$ &
        $\gh_{2n;\C} \oplus \gh_{2;\C}$ \\ \hline
19b & $U(n) \times U(2)$ &
        $\C^{n\times 2}\oplus\C^2$ & $\R\oplus\R$ & $0$ &
        $\gh_{2n;\C} \oplus \gh_{2;\C}$ \\ \hline
19c & $U(1)Sp(\tfrac{n}{2}) \times U(2)$ &
        $\C^{n\times 2}\oplus\C^2$ & $\R\oplus\R$ & $0$ &
        $\gh_{2n;\C} \oplus \gh_{2;\C}$ \\ \hline
20aa & $SU(n) \times SU(2)\times SU(m)$ &
        $\C^{n\times 2}\oplus\C^{2\times m}$ &
        $\R \oplus \R$ & $0$ & $\gh_{2n;\C} \oplus \gh_{2m;\C}$ \\ \hline
20ab & $SU(n) \times SU(2)\times U(m)$ &
        $\C^{n\times 2}\oplus\C^{2\times m}$ &
        $\R \oplus \R$ & $0$ & $\gh_{2n;\C} \oplus \gh_{2m;\C}$ \\ \hline
20ac & $SU(n) \times SU(2)\times U(1)Sp(\tfrac{m}{2})$ &
        $\C^{n\times 2}\oplus\C^{2\times m}$ &
        $\R \oplus \R$ & $0$ & $\gh_{2n;\C} \oplus \gh_{2m;\C}$ \\ \hline
20ba & $U(n) \times SU(2)\times SU(m)$ &
        $\C^{n\times 2}\oplus\C^{2\times m}$ &
        $\R \oplus \R$ & $0$ & $\gh_{2n;\C} \oplus \gh_{2m;\C}$ \\ \hline
20bb & $U(n) \times SU(2)\times U(m)$ &
        $\C^{n\times 2}\oplus\C^{2\times m}$ &
        $\R \oplus \R$ & $0$ & $\gh_{2n;\C} \oplus \gh_{2m;\C}$ \\ \hline
20bc & $U(n) \times SU(2)\times U(1)Sp(\tfrac{m}{2})$ &
        $\C^{n\times 2}\oplus\C^{2\times m}$ &
        $\R \oplus \R$ & $0$ & $\gh_{2n;\C} \oplus \gh_{2m;\C}$ \\ \hline
20ca & $U(1)Sp(\tfrac{n}{2}) \times SU(2)\times SU(m)$ &
        $\C^{n\times 2}\oplus\C^{2\times m}$ &
        $\R \oplus \R$ & $0$ & $\gh_{2n;\C} \oplus \gh_{2m;\C}$ \\ \hline
20cb & $U(1)Sp(\tfrac{n}{2}) \times SU(2)\times U(m)$ &
        $\C^{n\times 2}\oplus\C^{2\times m}$ &
        $\R \oplus \R$ & $0$ & $\gh_{2n;\C} \oplus \gh_{2m;\C}$ \\ \hline
20cc & $U(1)Sp(\tfrac{n}{2}) \times SU(2)\times$ &
        $\C^{n\times 2}\oplus\C^{2\times m}$ &
        $\R \oplus \R$ & $0$ & $\gh_{2n;\C} \oplus \gh_{2m;\C}$ \\
   & \phantom{$\{SU(m),U(m)\}$} $U(1)Sp(\tfrac{m}{2})\}$ & & & & \\ \hline
21a & $SU(n) \times SU(2) \times U(4)$ &
        $\C^{n\times 2}\oplus\C^{2\times 4}$ &
        $\R \oplus \R$ & $\R^6$ & 
	$\gh_{2n;\C} \oplus \gh_{8;\C}$ \\ \hline
21b & $U(n) \times SU(2) \times U(4)$ &
        $\C^{n\times 2}\oplus\C^{2\times 4}$ &
        $\R \oplus \R$ & $\R^6$ & 
	$\gh_{2n;\C} \oplus \gh_{8;\C}$ \\ \hline
21c & $U(1)Sp(\tfrac{n}{2}) \times SU(2) \times U(4)$ &
        $\C^{n\times 2}\oplus\C^{2\times 4}$ &
        $\R \oplus \R$ & $\R^6$ & 
	$\gh_{2n;\C} \oplus \gh_{8;\C}$ \\ \hline
\end{tabular}
\end{equation}
}

By inspection of each row of the table we arrive at
\begin{proposition}\label{yak1}
Each of the $28$ strict 
direct systems $\{(G'_\ell,K_\ell)\}$ 
of {\rm Table \ref{ind-vin-table-ipms}} satisfies {\rm (\ref{nil-conditions})}.
\end{proposition}

If $\ell \leqq k$ in the index set then $\gz'_\ell \hookrightarrow \gz'_k$ is
surjective.  We identify each of the $\gz'_\ell$ with 
$\gz' := \varinjlim \gz'_\ell$.  As in Lemma \ref{v-basis}
we have $\ga'_\ell$, the zero set of the polynomial $\Pf(b_{\ell,t'})$ on
$\gz'_\ell$, and $\ga' := \bigcup \ga'_\ell$ is a set of measure zero
in $(\gz')^*$.  And as in (\ref{T}) we denote 
$T' = \{t' \in (\gz')^* \mid \text{ each } \Pf(b_{\ell,t'}) \ne 0\}
= (\gz')^* \setminus \ga'$.  
\medskip

If $t' \in T'$ then $K_{\ell,t'}$, $G'_{\ell,t'} = N_\ell \rtimes K'_{\ell,t'}$
and $G_{\ell,t'} = N_\ell \rtimes K_{\ell,t'}$ are its respective stabilizers 
in $K_\ell$, $G'_\ell$ and $G_\ell$.  Theorem \ref{iso-gelfand} tells us that
$(G'_{\ell,t'}, K_{\ell,t'})$ is a Gelfand pair, so in particular the
action of $K_{\ell,t'}$ on $\C[\gv_\ell]$ is multiplicity free.  As
$G_{\ell,t'}$ and $G'_{\ell,t'}$ correspond to the same multiplicity free
action of $K_{\ell,t'}$ on $\C[\gv_\ell]$, Carcano's Theorem says that
$(G_{\ell,t'}, K_{\ell,t'})$ is a Gelfand pair.  Now Lemma \ref{k-invariant}
and Proposition \ref{restrict-invariants} apply.  As above, this leads to

\begin{theorem}\label{big-hhnt}
Let $\{(G'_\ell,K_\ell)\}$ be one of the twenty eight direct systems of
{\rm Table \ref{ind-vin-table-ipms}}.  Denote $G' = \varinjlim G'_\ell$ and
$K = \varinjlim K_\ell$.  Then the unitary direct system 
$\{L^2(G'_\ell), \zeta'_{\ell,\widetilde{\ell}}\}$ given by that
of {\rm Theorem \ref{hnt3}}, restricts to a unitary direct system
$\{L^2(G'_\ell/K_\ell), \zeta'_{\ell,\widetilde{\ell}}\}$, the
Hilbert space $L^2(G'/K) := 
\varinjlim \{L^2(G'_\ell/K_\ell), \zeta'_{\ell,\widetilde{\ell}}\}$ 
is the subspace of
$L^2(G') := \varinjlim \{L^2(G'_\ell), \zeta'_{\ell,\widetilde{\ell}}\}$ 
consisting of right--$K$--invariant
functions, and the natural unitary representation of $G'$ on
$L^2(G'/K)$ is a multiplicity free direct integral of lim--irreducible
representations.
\end{theorem}

Now define $\gz'' = \varinjlim \gz''_\ell$ and let $Z'' := \varinjlim Z''_\ell$
denote the corresponding vector group.  In Table \ref{ind-vin-table-ipms}
the $\gz''_\ell$ are constant except for the entries of row 1, 
where $\gz'' = \gs\gu(\infty)$, and rows 11a,b,c, where $\gz''$ can be
to any of $\H^{n\times\infty}$, $\H^{\infty \times m}$ and 
$\H^{\infty\times\infty}$.  In any case, 
$(\gz'')^* = \varprojlim (\gz''_\ell)^*$.  
\medskip

In the cases where the $\gz''_\ell$ are zero, i.e. $\gz'' = 0$,
we have $G_\ell = G'_\ell$, so the natural unitary representation of 
$G = \varinjlim G_\ell$
on $L^2(G/K)$ is  multiplicity free by Theorem \ref{big-hhnt}.
\medskip

The cases where the $\gz''_\ell$ are nonzero but constant are
{\footnotesize
\begin{equation} \label{ind-vin-table-ipms-nzc}
\begin{tabular}{|c|l|l|l|l|l|}\hline
 & Group $K_\ell$ & $K_\ell$--module $\gv_\ell$ & module $\gz_\ell' = [\gn_\ell,\gn_\ell]$ &
        module $\gz''_\ell$ & Algebra $\gn'_\ell$
        \\ \hline
6 & $S(U(4)\times U(m))$ & $\C^{4 \times m}$ & $\R$ & $\R^6$ &
        $\gh_{4m;\C}$ \\ \hline
18a & $SU(n) \times SU(2)$ &
        $\C^{n\times 2}$ & $\R$ & $\gs\gu(2)$ &
        $\gh_{2n;\C}$ \\ \hline
18b & $U(n) \times SU(2)$ &
        $\C^{n\times 2}$ & $\R$ & $\gs\gu(2)$ &
        $\gh_{2n;\C}$ \\ \hline
18c & $U(1)Sp(\tfrac{n}{2}) \times SU(2)$ &
        $\C^{n\times 2}$ & $\R$ & $\gs\gu(2)$ &
        $\gh_{2n;\C}$ \\ \hline
21a & $SU(n) \times SU(2) \times U(4)$ &
        $\C^{n\times 2}\oplus\C^{2\times 4}$ &
        $\R \oplus \R$ & $\R^6$ &
        $\gh_{2n;\C} \oplus \gh_{8;\C}$ \\ \hline
21b & $U(n) \times SU(2) \times U(4)$ &
        $\C^{n\times 2}\oplus\C^{2\times 4}$ &
        $\R \oplus \R$ & $\R^6$ &
        $\gh_{2n;\C} \oplus \gh_{8;\C}$ \\ \hline
21c & $U(1)Sp(\tfrac{n}{2}) \times SU(2) \times U(4)$ &
        $\C^{n\times 2}\oplus\C^{2\times 4}$ &
        $\R \oplus \R$ & $\R^6$ &
        $\gh_{2n;\C} \oplus \gh_{8;\C}$ \\ \hline
\end{tabular}
\end{equation}
}
In those cases $\{(G_\ell,K_\ell)\}$ satisfies (\ref{nil-conditions}), so 
the considerations leading to Theorem \ref{big-hhnt} apply directly to the 
$\{(G_\ell,K_\ell)\}$.
We turn now to the other cases.  They are given by
{\footnotesize
\begin{equation} \label{ind-vin-table-ipms-ubd}
\begin{tabular}{|c|l|l|l|l|l|}\hline
 & Group $K_\ell$ & $K_\ell$--module $\gv_\ell$ & module $\gz_\ell' = [\gn_\ell,\gn_\ell]$ &
	module $\gz''_\ell$ & Algebra $\gn'_\ell$
        \\ \hline
1 & $U(n)$ & $\C^n$ & $\R$ & $\gs\gu(n)$ &
        $\gh_{n;\C}$  \\ \hline
11a & $Sp(n)\times Sp(1)\times Sp(m)$ & $\H^n$ & $\Im\H$ &
        $\H^{n\times m}$ & $\gh_{n;\H}$ \\ \hline
11b & $Sp(n)\times U(1)\times Sp(m)$ & $\H^n$ & $\Im\H$ &
        $\H^{n\times m}$ & $\gh_{n;\H}$ \\ \hline
11c & $Sp(n)\times \{1\} \times Sp(m)$ & $\H^n$ & $\Im\H$ &
        $\H^{n\times m}$ & $\gh_{n;\H}$ \\ \hline
\end{tabular}
\end{equation}
}

For entry 1 in Table \ref{ind-vin-table-ipms-ubd} we have 
$T' = (\gz')^*\setminus \{0\} = \R \setminus \{0\}$.  The generic orbits of 
$K_\ell = U(n)$ on $\gz''_\ell = \gs\gu(n)$ are those for which all
eigenvalues of the matrix in $\gs\gu(n)$ are distinct, so the
generic $t = (t',t'') \in T$ are those for which the stabilizer
$K_{\ell,t}$ is a maximal torus in $K_\ell$.  Here note that the irreducible
subspaces for $K_{\ell,t}$ on $\C[\gv_\ell] = \C[\C^n]$ each consists of
the multiples of a monomial, so the action of $K_{\ell,t}$ on $\C[\gv_\ell]$
is multiplicity free.  Now, exactly as in the considerations leading up
to Theorem \ref{big-hnt}, the natural unitary representation of $G$ on
$L^2(G/K)$ is a multiplicity free direct integral of lim--irreducible
representations.
\medskip

The argument is more or less the same for entries 11a,b,c in 
Table \ref{ind-vin-table-ipms-ubd}.  
Here $T' = (\gz')^*\setminus \{0\} = \Im\H \setminus \{0\}$ where $\Im\H$
is identified with its real dual using the inner product 
$\langle z'_1,z'_2\rangle = \Re (z'_1 \overline{z'_2})$.  
The action of
$K_\ell = Sp(n)\times \{Sp(1), U(1), \{1\}\}\times Sp(m)$ on
$\gv_\ell = \H^n$ is $(k_1,k_2,k_3):v \mapsto k_1v\overline{k_2}$, on
$\gz'_\ell = \Im\H$ is $(k_1,k_2,k_3): z' \mapsto k_2z'\overline{k_2}$, and on 
$\gz''_\ell = \H^{n\times m}$ is $(k_1,k_2,k_3): z'' \mapsto k_1z''k_3^*$.
Here the composition $\H^n\times\H^n \to \Im\H$ is $(v,w) \mapsto \Im v^*w$
and $k_2$ is a quaternionic scalar.  The stabilizer of $t' \in T'$ in 
$K_\ell$ is $Sp(n)\times \{U(1), U(1)\text{ or }\{\pm 1\}, \{1\}\}\times Sp(m)$;
generically that stabilizer is 
$Sp(n)\times \{U(1),\{\pm 1\}, \{1\}\}\times Sp(m)$.  Again, as in the 
considerations leading up to Proposition \ref{restrict-more-invariants} and 
Theorem \ref{big-hnt}, the natural unitary representation of $G$ on
$L^2(G/K)$ is a multiplicity free direct integral of lim--irreducible
representations.  In summary,

\begin{theorem}\label{big-hhhhnt}
Let $\{(G_\ell,K_\ell)\}$ be one of the twenty eight direct systems of
{\rm Table \ref{ind-vin-table-ipms}}.  Then the unitary direct system 
$\{L^2(G_\ell), \zeta_{\ell,\widetilde{\ell}}\}$ analogous to that of 
{\rm Theorem \ref{hnt3}} restricts to a unitary direct system 
$\{L^2(G_\ell/K_\ell), \zeta_{\ell,\widetilde{\ell}}\}$.  The
Hilbert space $L^2(G/K) := 
\varinjlim \{L^2(G_\ell/K_\ell), \zeta_{\ell,\widetilde{\ell}}\}$ 
is the subspace of $L^2(G) := 
\varinjlim \{L^2(G_\ell), \zeta_{\ell,\widetilde{\ell}}\}$ consisting of 
right--$K$--invariant functions, and the natural unitary representation of 
$G$ on $L^2(G/K)$ is a multiplicity free direct integral of lim--irreducible
representations.
\end{theorem}

\section*{Appendix A: Formal Degrees of Induced Representations}
\label{appendixa}
\setcounter{equation}{0}

In this section we work out the formal degree of an irreducible induced
representation $\Ind_M^L(\gamma)$, where $\gamma$ is a square integrable
representation of $L$, $L/M$ is compact, and $L/M$ has a positive
$L$--invariant measure.  The result, which is
suggested by (\ref{formaldeg}), is not surprising, but does not seem to
be in the literature.
\medskip

Let $L$ be a separable locally compact group, $M$ a closed
subgroup of $L$, and $J$ a closed central subgroup 
of $M$ that is normal in $L$.  Suppose that $\gamma$ is an irreducible 
square integrable (modulo $J$) unitary representation of $M$.  Then 
$\gamma$ has well defined formal degree $\deg\gamma$ in
the usual sense:  if $u$, $v$, $u'$ and $v'$ belong to the representation
space $\cH_\gamma$, and if we write $f_{u,v}$ for the 
coefficient $f_{u,v}(m) = \langle u, \gamma(m)v\rangle_{\cH_\gamma}$, then
$\langle f_{u,v}, f_{u',v'} \rangle_{L^2(M/J)} = \tfrac{1}{\deg \gamma}
	\langle u, u'\rangle \overline{\langle v, v'\rangle}$.
\medskip

The modular functions $\Delta_M$ and $\Delta_L$ coincide on $M$.  This
is just another way of saying that
we have a unique (up to scale) $L$--invariant Radon measure $d(\ell M)$
on $L/M$, and thus an $L$--invariant integral 
$f \mapsto \int_{L/M} f(\ell)d(\ell M)$ for functions $f \in C_c(L/M)$.
\medskip

Denote $\widetilde{\gamma} = \Ind_M^L(\gamma)$ and let 
$\widetilde{\cH_\gamma}$ be its representation space.  The elements of 
$\widetilde{\cH_\gamma}$ are the measurable functions 
$\varphi: L \to \cH_\gamma$ such that (i) 
$\varphi(\ell m) = \gamma(m)^{-1}\varphi(\ell)$ (for $\ell \in L$
and $m \in M$) and (ii) $\int_{L/M} ||\varphi(\ell)||^2d(\ell M) < \infty$.  
The action $\widetilde{\gamma}$ of $L$ on $\widetilde{\cH_\gamma}$ is
$[\widetilde{\gamma}(\ell)\varphi](\ell') = \varphi(\ell^{-1}\ell')$. 
The inner product on $\widetilde{\cH_\gamma}$ is
$\langle \varphi, \psi \rangle_{_{\widetilde{\cH_\gamma}}} = 
\int_{L/M} \langle \varphi(\ell), \psi(\ell)\rangle_{_{\cH_\gamma}} d(\ell M)$.
\medskip

The group $M$ acts on the fiber $\cH_\gamma$ of 
$\widetilde{\cH_\gamma} \to L/M$ at $1 M$ by $\gamma$, and more generally
the stabilizer $\ell M \ell^{-1}$ of the fiber $\ell \cH_\gamma$ at $\ell M$ 
acts on that fiber by 
$\gamma_\ell (\ell m \ell^{-1})(\ell v) = \ell \cdot\gamma(m)v$.  Since $J$ is
normal in $L$ it sits in $\ell M \ell^{-1}$ and $\gamma_\ell|_J$ consists
of scalar transformations of $\ell \cH_\gamma$.  Now the representations
$\gamma_\ell$ all are square integrable mod $J$ and have the same formal 
degree $\deg \gamma$.
\medskip

If $u,v \in \cH_\gamma$ we have the coefficient function
$f_{u,v}(m) = \langle u, \gamma(m)v \rangle_{_{\cH_\gamma}}$ on $M$.  If
$\varphi, \psi \in \widetilde{\cH_\gamma}$ we have the coefficient
$\widetilde{f}_{\varphi, \psi}(\ell) 
= \langle \varphi, \widetilde{\gamma}(\ell) 
\psi \rangle_{_{\widetilde{\cH_\gamma}}}$ on $L$.  Now let 
$\varphi, \psi \in \widetilde{\cH_\gamma}$ such that, as functions
from $L$ to $\cH_\gamma$, both $\varphi$ and $\psi$ are continuous.  Their 
support is compact modulo $M$ because $L/M$ is compact.  Now 
$\widetilde{f}_{\varphi, \psi}:L \to \C$
has support that is compact modulo $M$.  For every $\ell \in L$, 
$\widetilde{f}_{\varphi, \psi}|_{\ell M}$ is a matrix coefficient of 
$\gamma_\ell$ and $\left |\widetilde{f}_{\varphi, \psi}|_{\ell M}\right |^2$ 
has integral (integrate over $\ell M \ell^{-1}/J$) equal to 
$\frac{1}{\deg\gamma}||\varphi(\ell)||^2_{_{\cH_\gamma}}
	||\psi(\ell)||^2_{_{\cH_\gamma}}$.
The functions $\ell \mapsto ||\varphi(\ell)||^2_{_{\cH_\gamma}}$
and $\ell \mapsto ||\psi(\ell)||^2_{_{\cH_\gamma}}$ are continuous on $L/M$, so 
their $L^2(L/M)$ inner product converges.  Now
$$
\begin{aligned}
\int_{L/J} |\widetilde{f}_{\varphi, \psi} (\ell)|^2 d(\ell M) &= 
\int_{L/M}\left ( \int_{\ell M}|\widetilde{f}_{\varphi, \psi}(\ell m)|^2
	d(mJ) \right ) d(\ell M) \\
&= \frac{1}{\deg\gamma} \int_{L/M} ||\varphi(\ell)||^2_{_{\cH_\gamma}}
	||\psi(\ell)||^2_{_{\cH_\gamma}} d(\ell M) < \infty.
\end{aligned}
$$
Thus $\widetilde{\gamma}$ has a nonzero square integrable (mod $J$) coefficient.
Since it is irreducible, all its coefficients are square integrable 
(mod $J$).  In particular the formal degree $\deg\widetilde{\gamma}$ is
defined.  In summary,
\medskip

\noindent {\bf Theorem A.1\ \ }{\em
Let $L$ be a separable locally compact group, $M \subset L$ a closed
subgroup, and $J$ a central subgroup of $M$ that is normal in $L$.  Suppose
that $L/M$ is compact and has a nonzero $L$--invariant Radon measure.
Let $\gamma$ be a square integrable {\rm (modulo $J$)} irreducible unitary 
representation of $M$ and $\deg \gamma$ its formal degree such
that $\widetilde{\gamma}:= \Ind_M^L(\gamma)$ is irreducible.  Then
$\widetilde{\gamma}$ is square integrable {\rm (modulo $J$)}
and it has a well defined formal degree $\deg\widetilde{\gamma}$.
}
\medskip

Of course, if $L/M$ is finite, then Theorem A.1 becomes
trivial, and there if we use counting measure on $L/M$ then 
$\deg\widetilde{\gamma} = |L/M|\deg\gamma$.  However, in effect we use
the result in (\ref{formaldeg}) where $L/M$ is compact but infinite.

\section*{Appendix B: A Computational Argument For Theorem \ref{iso-gelfand}}
\label{appendixb}

In this appendix we give a computational proof of Theorem \ref{iso-gelfand} 
for the spaces of Table \ref{ind-vin-table}.  This provides somewhat more
information and could be useful in studying the spherical functions. 
In order to align the $K_n$--invariants in $L^2(G_n)$ and pass to 
$L^2(G_n/K_n)$ we will need

\noindent {\bf Theorem B.1\ \ }{\em
Let $\{(G_n,K_n)\}$ be one of the thirteen direct systems 
{\rm Table \ref{ind-vin-table}} and let $t \in T$.  Then the representation
of $K_{n,t}$ on $\C[\gv_n]$ is multiplicity free.
}
\medskip

\noindent {\bf Proof.}  If $\dim \gz = 1$ the assertion follows from
Carcano's Theorem.  That leaves table entries 17, 18, 20a, 20b and 22.
In those cases we explicitly decompose $\C[\gv_n]$ under the $K_{n,t}$.
\medskip

In the case of entry 17, the big factor $K_n'' = Sp(n)$ is irreducible on
$\gv_n = \H^n = \C^{2n}$, and $K_{n,t} = U(1)\times Sp(n)$.  The
representation of $K_{n,t}$ on $\C[gv_n]$ is the same as that of the
Gelfand pair listed on row 4 of Table \ref{kac-table} and also row 4 of
Table \ref{jaw-table}.  It follows that the representation of $K_{n,t}$
on $\C[\gv_n]$ is multiplicity free.
\medskip

In the case of entry 20a we have $K_{n,t} = K'_{n,t} \times  SU(n)$ where
$K'_{n,t}$ is the circle group consisting of all 
$k_a = \left ( \begin{smallmatrix} a & 0 \\ 0 & 1/a \end{smallmatrix} \right )$ 
with $|a| = 1$.  Then $\gv_n = \C^n_+ \oplus
\C^n_-$ where $k_a$ acts on $\C^n_\pm$ as multiplication by $a^{\pm 1}$.
Define (i) $\P_{n,1,m_1}$ is the space of 
polynomials of degree $m_1$ on $\C^n_+$, (ii) $\P_{n,1,m_2}$ is the space of
polynomials of degree $m_2$ on $\C^n_-$, and (iii) $\P_{n,m} =
\P_{n,1,m_1} \otimes \P_{n,1,m_2}$ where $m = (m_1,m_2)$.  Then
$K_{n,t}$ acts on $\P_{n,1,m_1}$ by 
\setlength{\unitlength}{.65 mm}
\begin{picture}(80,10)
\put(20,1){\circle{2}}
\put(17,4){$m_1$}
\put(21,1){\line(1,0){13}}
\put(35,1){\circle{2}}
\put(36,1){\line(1,0){13}}
\put(52,1){\circle*{1}}
\put(55,1){\circle*{1}}
\put(58,1){\circle*{1}}
\put(61,1){\line(1,0){13}}
\put(75,1){\circle{2}}
\put(3,0){$\times$}
\put(2,4){$m_1$}
\end{picture}
and on $\P_{n,2,m_2}$ by
\setlength{\unitlength}{.65 mm}
\begin{picture}(80,10)
\put(20,1){\circle{2}}
\put(17,4){$m_2$}
\put(21,1){\line(1,0){13}}
\put(35,1){\circle{2}}
\put(36,1){\line(1,0){13}}
\put(52,1){\circle*{1}}
\put(55,1){\circle*{1}}
\put(58,1){\circle*{1}}
\put(61,1){\line(1,0){13}}
\put(75,1){\circle{2}}
\put(3,0){$\times$}
\put(0,4){$-m_2$}
\end{picture}.
Thus the representation of $K_{n,t}$ on $\P_{n,m}$ is
$\sum_{p+2q = m_1+m_2}\Bigl (
\setlength{\unitlength}{.65 mm}
\begin{picture}(105,10)
\put(30,1){\circle{2}}
\put(29,4){$p$} 
\put(31,1){\line(1,0){13}} 
\put(45,1){\circle{2}}
\put(44,4){$q$}
\put(46,1){\line(1,0){13}} 
\put(60,1){\circle{2}}
\put(61,1){\line(1,0){13}} 
\put(77,1){\circle*{1}} 
\put(80,1){\circle*{1}}
\put(83,1){\circle*{1}}
\put(86,1){\line(1,0){13}}
\put(100,1){\circle{2}}
\put(9,0){$\times$}
\put(0,4){$m_1-m_2$}
\end{picture}\Bigr )$.
Multiplicities in the representation of $K_{n,t}$ on $\C[\gv_n]$
would lead to equations $m'_1 - m'_2 = m_1 - m_2$, $p' = p$ and $q' = q$
with $p' + 2q' = m'_1 + m'_2$ and $p + 2q = m_1 + m_2$, forcing $m'_1 = m_1$
and $m'_2 = m_2$.  That is a contradiction because 
representation of $K_{n,t}$ on $\P_{n,m}$ is multiplicity free.  We
conclude that the representation of $K_{n,t}$ on $\C[\gv_n]$ is
multiplicity free.
\medskip

In the case of entry 20b, the representation of $K_{n,t}$ on $\C[\gv_n]$ is
multiplicity free as a consequence of the multiplicity free result for 
entry 20a.
\medskip

In the case of entry 22 we view $\gv_n$ as $\C^2 \otimes \C^{2n}$.  Here
$K_{n,t} = K'_{n,t} \times  Sp(n)$ where $K'_{n,t}$ is the $2$--torus 
group consisting of all
$k_{a,b} = \left (\begin{smallmatrix} a & 0 \\ 0 & b \end{smallmatrix}\right )$
with $|a| = |b| = 1$.  Then $\gv_n = \C^n_+ \oplus
\C^n_-$ where $k_{a,b}$ acts on $\C^n_+$ as multiplication by $a$ and on
$\C^n_-$ as multiplication by $b$.  The representation of $Sp(n)$ on either
of $\C^n_\pm$ is 
\setlength{\unitlength}{.60 mm}
\begin{picture}(80,7)
\put(5,1){\circle{2}}
\put(4,4){{\footnotesize $1$}}
\put(6,1){\line(1,0){13}}
\put(20,1){\circle{2}}
\put(21,1){\line(1,0){13}}
\put(37,1){\circle*{1}}
\put(40,1){\circle*{1}}
\put(43,1){\circle*{1}}
\put(46,1){\line(1,0){13}}
\put(60,1){\circle{2}}
\put(61,0.5){\line(1,0){13}}
\put(61,1.5){\line(1,0){13}}
\put(75,1){\circle{2}}
\put(66,-0.2){$<$}
\end{picture},
so the representation on polynomials of degree $d$ is the symmetric power
$S^d(\setlength{\unitlength}{.60 mm}
\begin{picture}(80,7)
\put(5,1){\circle{2}}
\put(4,4){{\footnotesize $1$}}
\put(6,1){\line(1,0){13}}
\put(20,1){\circle{2}}
\put(21,1){\line(1,0){13}}
\put(37,1){\circle*{1}}
\put(40,1){\circle*{1}}
\put(43,1){\circle*{1}}
\put(46,1){\line(1,0){13}}
\put(60,1){\circle{2}}
\put(61,0.5){\line(1,0){13}}
\put(61,1.5){\line(1,0){13}}
\put(75,1){\circle{2}}
\put(66,-0.2){$<$}
\end{picture}) =
\setlength{\unitlength}{.60 mm}
\begin{picture}(80,7)
\put(5,1){\circle{2}}
\put(4,4){{\footnotesize $d$}}
\put(6,1){\line(1,0){13}}
\put(20,1){\circle{2}}
\put(21,1){\line(1,0){13}}
\put(37,1){\circle*{1}}
\put(40,1){\circle*{1}}
\put(43,1){\circle*{1}}
\put(46,1){\line(1,0){13}}
\put(60,1){\circle{2}}
\put(61,0.5){\line(1,0){13}}
\put(61,1.5){\line(1,0){13}}
\put(75,1){\circle{2}}
\put(66,-0.2){$<$}
\end{picture}$.
Now the branching rule (with $r \leqq s$)
$$
\begin{aligned}
&\setlength{\unitlength}{.60 mm}
\begin{picture}(80,7)
\put(5,1){\circle{2}}
\put(4,4){{\footnotesize $r$}}
\put(6,1){\line(1,0){13}}
\put(20,1){\circle{2}}
\put(21,1){\line(1,0){13}}
\put(37,1){\circle*{1}}
\put(40,1){\circle*{1}}
\put(43,1){\circle*{1}}
\put(46,1){\line(1,0){13}}
\put(60,1){\circle{2}}
\put(61,0.5){\line(1,0){13}}
\put(61,1.5){\line(1,0){13}}
\put(75,1){\circle{2}}
\put(66,-0.2){$<$}
\end{picture}
\otimes
\setlength{\unitlength}{.60 mm}
\begin{picture}(80,7)
\put(5,1){\circle{2}}
\put(4,4){{\footnotesize $s$}}
\put(6,1){\line(1,0){13}}
\put(20,1){\circle{2}}
\put(21,1){\line(1,0){13}}
\put(37,1){\circle*{1}}
\put(40,1){\circle*{1}}
\put(43,1){\circle*{1}}
\put(46,1){\line(1,0){13}}
\put(60,1){\circle{2}}
\put(61,0.5){\line(1,0){13}}
\put(61,1.5){\line(1,0){13}}
\put(75,1){\circle{2}}
\put(66,-0.2){$<$}
\end{picture}\\
&\phantom{XXXXXX} = 
\sum_{v=0}^{r-s} \ \ \sum_{u=0}^v \phantom{X} (
\setlength{\unitlength}{.60 mm}
\begin{picture}(105,7)
\put(15,1){\circle{2}}
\put(4,4){{\tiny $r+s-2v$}}
\put(16,1){\line(1,0){13}}
\put(30,1){\circle{2}}
\put(29,4){{\tiny $u$}}
\put(31,1){\line(1,0){13}}
\put(45,1){\circle{2}}
\put(46,1){\line(1,0){13}}
\put(62,1){\circle*{1}}
\put(65,1){\circle*{1}}
\put(68,1){\circle*{1}}
\put(71,1){\line(1,0){13}}
\put(85,1){\circle{2}}
\put(86,0.5){\line(1,0){13}}
\put(86,1.5){\line(1,0){13}}
\put(100,1){\circle{2}}
\put(91,-0.2){$<$}
\end{picture} )
\end{aligned}
$$
shows that the representation of $K_{n,t}$ on $\P_{n,m}$ is
$$
\sum_{v=0}^{|m_1-m_2|} \ \ \sum_{u=0}^v \phantom{X} (
\setlength{\unitlength}{.60 mm}
\begin{picture}(145,7)
\put(5,-1){$\times$}
\put(5,4){{\tiny $m_1$}}
\put(20,-1){$\times$}
\put(20,4){{\tiny $m_2$}}
\put(50,1){\circle{2}}
\put(35,4){{\tiny $m_1+m_2-2v$}}
\put(51,1){\line(1,0){18}}
\put(70,1){\circle{2}}
\put(69,4){{\tiny $u$}}
\put(71,1){\line(1,0){13}}
\put(85,1){\circle{2}}
\put(86,1){\line(1,0){13}}
\put(102,1){\circle*{1}}
\put(105,1){\circle*{1}}
\put(108,1){\circle*{1}}
\put(111,1){\line(1,0){13}}
\put(125,1){\circle{2}}
\put(126,0.5){\line(1,0){13}}
\put(126,1.5){\line(1,0){13}}
\put(140,1){\circle{2}}
\put(131,-0.2){$<$}
\end{picture} ).
$$
If two such irreducible summands are equivalent, say for 
$(m_1,m_2,u,v)$ and $(m'_1,m'_2,u',v')$, then evidently $m_1 = m'_1$,
$m_2 = m'_2$, $m_1+m_2-2v = m'_1 + m'_2 - 2v'$ and $u = u'$.  We conclude
that the representation of $K_{n,t}$ on $\C[\gv_n]$ is multiplicity free.
\medskip

In the case of entry 18 the $K_{n,t}$ correspond to the centralizers of tori
(of dimensions $0$, $1$ or $2$) in $Sp(2)$.  If we view $Sp(2)$ from its diagram
\setlength{\unitlength}{.60 mm}
\begin{picture}(25,7)
\put(5,1){\circle{2}}
\put(4,4){{\footnotesize $\alpha$}}
\put(6,0.5){\line(1,0){13}}
\put(6,1.5){\line(1,0){13}}
\put(20,1){\circle{2}}
\put(19,4){{\footnotesize $\beta$}}
\put(11,-0.2){$<$}
\end{picture}
the centralizers $K_{n,t} = K'_{n,t} \times Sp(n)$ are given up to conjugacy 
by $K'_{n,t} = Sp(2)$, $K'_{n,t} \cong U(2)$ with simple root $\alpha$, 
$K'_{n,t} \cong U(2)$ with simple root $\beta$, and $K'_{n,t} = T^2$ maximal
torus of $Sp(2)$.  Specifically, if we realize $\gs\gp(2)$ as the space of
$2\times 2$ quaternionic matrices $\xi$ with $\xi+\xi^* = 0$, it has basis
$$
(\begin{smallmatrix}i & 0 \\ 0 & 0 \end{smallmatrix}), \ 
 (\begin{smallmatrix}j & 0 \\ 0 & 0 \end{smallmatrix}), \ 
 (\begin{smallmatrix}k & 0 \\ 0 & 0 \end{smallmatrix}), \ 
(\begin{smallmatrix}0 & 0 \\ 0 & i \end{smallmatrix}), \ 
 (\begin{smallmatrix}0 & 0 \\ 0 & j \end{smallmatrix}), \ 
 (\begin{smallmatrix}0 & 0 \\ 0 & k \end{smallmatrix}), \ 
(\begin{smallmatrix}0 & 1 \\ -1 & 0 \end{smallmatrix}), \
 (\begin{smallmatrix}0 & i \\ i & 0 \end{smallmatrix}), \
 (\begin{smallmatrix}0 & j \\ j & 0 \end{smallmatrix}) \text{ and }
 (\begin{smallmatrix}0 & k \\ k & 0 \end{smallmatrix}). 
$$
Write the complexifying $\C$ as $\R + b\R$ where $b^2 = -1$, to avoid 
the notation of the quaternions used to express $\gs\gp(2)$, and let
$\varepsilon_1, \varepsilon_2$ be the usual linear functional (that
pick out the diagonal entries of a matrix).  Using the Cartan subalgebra
spanned by $(\begin{smallmatrix}i & 0 \\ 0 & 0 \end{smallmatrix})$ and
$(\begin{smallmatrix}0 & 0 \\ 0 & i \end{smallmatrix})$,
$$
\begin{aligned}
&(\begin{smallmatrix}j \pm bk & 0 \\ 0 & 0 \end{smallmatrix})
	\text{ is a root vector for } \pm\varepsilon_1 \text{ and }
(\begin{smallmatrix} 0 & 0 \\ 0 & j \pm bk \end{smallmatrix})
	\text{ is a root vector for } \pm\varepsilon_2; \text{ and } \\
\bigskip
&(\begin{smallmatrix}0 & 1 \pm bi \\ -1 \pm bi & 0 \end{smallmatrix})
	\text{ is a root vector for } 
	\pm (\varepsilon_1 - \varepsilon_2) \text{ and }
(\begin{smallmatrix}0 & j \pm bk \\ j \pm bk & 0 \end{smallmatrix})
	\text{ is a root vector for } 
	\pm (\varepsilon_1 + \varepsilon_2).
\end{aligned}
$$
Thus the centralizer of 
$z_a := (\begin{smallmatrix}a_1 & 0 \\ 0 & a_2 \end{smallmatrix})$
is the direct sum of the Cartan subalgebra with the span of

(if $a_1 = 0 \ne a_2$): 
 $(\begin{smallmatrix}j & 0 \\ 0 & 0 \end{smallmatrix})$ and
 $(\begin{smallmatrix}k & 0 \\ 0 & 0 \end{smallmatrix})$, \ so
the centralizer of $z_a$ is 
 $(\begin{smallmatrix}\gs\gp(1) & 0 \\ 0 & 0 \end{smallmatrix})\oplus 
  (\begin{smallmatrix}0 & 0 \\ 0 & i\R \end{smallmatrix})$;

(if $a_1 \ne 0 = a_2$): 
 $(\begin{smallmatrix}0 & 0 \\ 0 & j \end{smallmatrix})$ and 
 $(\begin{smallmatrix}0 & 0 \\ 0 & k \end{smallmatrix})$,  \ so
the centralizer of $z_a$ is 
 $(\begin{smallmatrix} 0 & 0 \\ 0 & \gs\gp(1) \end{smallmatrix})\oplus 
  (\begin{smallmatrix} i\R & 0 \\ 0 & 0 \end{smallmatrix})$;

(if $0 \ne a_1 \ne a_2 \ne 0$):
 the centralizer of $z_a$ is the Cartan subalgebra, span of
 $(\begin{smallmatrix}i & 0 \\ 0 & 0 \end{smallmatrix})$ and
 $(\begin{smallmatrix}0 & 0 \\ 0 & i \end{smallmatrix})$;

(if $a_1 = a_2 \ne 0$):
 $(\begin{smallmatrix}0 & 1 \\ -1 & 0 \end{smallmatrix})$ and
 $(\begin{smallmatrix}0 & i \\ i & 0 \end{smallmatrix})$,  \ so
the centralizer of $z_a$ is $\gu(2)$ as $2 \times 2$ complex matrices;

(if $a_1 = -a_2 \ne 0$):
 $(\begin{smallmatrix}0 & j \\ j & 0 \end{smallmatrix})$ and
 $(\begin{smallmatrix}0 & k \\ k & 0 \end{smallmatrix})$,  \ so
the centralizer of $z_a$ is isomorphic to $\gu(2)$.
\medskip

\noindent Here the first two cases are $Sp(2)$--conjugate, and the last
two cases are $Sp(2)$--conjugate, because the Weyl group consists of all
signed permutations of the $\epsilon_\ell$.  Thus we really only have three
cases up to $Sp(2)$--conjugacy.  The nilpotent group $N_n$ is the (very)
generalized Heisenberg group $H_{2,n,0;\H}$ of \cite[Section 14.3]{W2},
so \cite[Theorem 14.3.1]{W2} shows that $z_a$ corresponds to a square
integrable representation of $N_n$ precisely when $a_1 \ne 0 \ne a_2$.
Now we need only consider the cases 
(i) $a_1 = a_2 \ne 0$ and (ii) $0 \ne a_1 \ne a_2 \ne 0$.
\medskip

Suppose $a_1 = a_2 \ne 0$.  Then the corresponding 
$K_{n,t} = SU(2)\times Sp(n)$.  Since $\dim_{_\R}\gv_n = 
\dim_{_\R} \H^{2\times n} = 8n$ we are now viewing $\gv$ as
$\C^{2 \times 2n}$, and the representation of $K_{n,t}$ (as a
subgroup of $U(n)$) on $\gv_n$ is the irreducible 
$(k',k''): v \mapsto k'v(k'')^{-1}$.  As $(N_n\rtimes K_{n,t}, K_{n,t})$
is a Gelfand pair, Carcano's Theorem says that 
the representation of $K_{n,t}$ on $\C[\gv_n]$ is multiplicity free.
\medskip

Suppose $0 \ne a_1 \ne a_2 \ne 0$.  Then $K_{n,t} = U(1) \times U(1)
\times Sp(n)$ and we view $\gv_n$ as $\C^{2\times 2n} = 
\C_1^{2n} \oplus \C_2^{2n}$ where the action $\tau$ of $K_{n,t}$ is
$\tau(k'_1,k'_2,k'') (v_1 \oplus v_2) = (k'_1v_1(k'')^{-1} \oplus
k'_2v_2(k'')^{-1})$.  The representation $\tau$ has diagram
$
(\setlength{\unitlength}{.60 mm}
\begin{picture}(90,7)
\put(4,-1){$\times$}
\put(6,4){{\tiny $1$}}
\put(13,-1){$\times$}
\put(25,1){\circle{2}}
\put(24,4){{\tiny $1$}}
\put(26,1){\line(1,0){8}}
\put(35,1){\circle{2}}
\put(36,1){\line(1,0){8}}
\put(45,1){\circle{2}}
\put(46,1){\line(1,0){8}}
\put(57,1){\circle*{1}}
\put(60,1){\circle*{1}}
\put(63,1){\circle*{1}}
\put(66,1){\line(1,0){8}}
\put(75,1){\circle{2}}
\put(76,0.5){\line(1,0){8}}
\put(76,1.5){\line(1,0){8}}
\put(85,1){\circle{2}}
\put(78,-0.2){$<$}
\end{picture})
\oplus
(\setlength{\unitlength}{.60 mm}
\begin{picture}(90,7)
\put(4,-1){$\times$}
\put(13,-1){$\times$}
\put(15,4){{\tiny $1$}}
\put(25,1){\circle{2}}
\put(24,4){{\tiny $1$}}
\put(26,1){\line(1,0){8}}
\put(35,1){\circle{2}}
\put(36,1){\line(1,0){8}}
\put(45,1){\circle{2}}
\put(46,1){\line(1,0){8}}
\put(57,1){\circle*{1}}
\put(60,1){\circle*{1}}
\put(63,1){\circle*{1}}
\put(66,1){\line(1,0){8}}
\put(75,1){\circle{2}}
\put(76,0.5){\line(1,0){8}}
\put(76,1.5){\line(1,0){8}}
\put(85,1){\circle{2}}
\put(78,-0.2){$<$}
\end{picture}).
$
Now decomposition of the symmetric powers $S^m(\tau)$
goes exactly as in the case of entry 22, and we conclude that
the representation of $K_{n,t}$ on $\C[\gv_n]$ is multiplicity free.  
\medskip

We have shown, by direct computation for each of the series of 
Table \ref{ind-vin-table}, that
the representation of $K_{n,t}$ on $\C[\gv_n]$ is multiplicity free.
\hfill $\square$
\medskip

One can carry out similar calculations for the series of Table 
\ref{ind-vin-table-ipms} showing by direct computation that in
each case the representation of $K_{n,t}$ on $\C[\gv_n]$ is multiplicity free.
\vfill\pagebreak

\end{document}